
\input amssym.def
\input amssym.tex


\def\item#1{\vskip1.3pt\hang\textindent {\rm #1}}


\tolerance=300
\pretolerance=200
\hfuzz=1pt
\vfuzz=1pt


\hoffset=0.6in
\voffset=0.8in

\hsize=5.8 true in 


\vsize=8.5 true in
\parindent=25pt
\mathsurround=1pt
\parskip=1pt plus .25pt minus .25pt
\normallineskiplimit=.99pt

\countdef\revised=100
\mathchardef\emptyset="001F 
\chardef\ss="19
\def\3{\ss}
\def\anf{$\lower1.2ex\hbox{"}$}
\def\frac#1#2{{#1 \over #2}}
\def\>{>\!\!>}
\def\<{<\!\!<}

\def\into{\hookrightarrow}
\def\ssarr{\hbox to 30pt{\rightarrowfill}}
\def\sarr{\hbox to 40pt{\rightarrowfill}}
\def\arr{\hbox to 60pt{\rightarrowfill}}
\def\larr{\hbox to 60pt{\leftarrowfill}}
\def\Arr{\hbox to 80pt{\rightarrowfill}}

{}

\def\ad{\mathop{\rm ad}\nolimits}

\def\Ad{\mathop{\rm Ad}\nolimits}

\def\Aut{\mathop{\rm Aut}\nolimits}

\def\cone{\mathop{\rm cone}\nolimits}
\def\conv{\mathop{\rm conv}\nolimits}

\def\Exp{\mathop{\rm Exp}\nolimits}

\def\Gl{\mathop{\rm Gl}\nolimits}

\def\HSp{\mathop{\rm HSp}\nolimits}
\def\Hol{\mathop{\rm Hol}\nolimits}%
%

\def\Im{\mathop{\rm Im}\nolimits}

\def\Inn{\mathop{\rm Inn}\nolimits}
\def\Int{\mathop{\rm int}\nolimits}

\def\Re{\mathop{\rm Re}\nolimits}

\def\Sl{\mathop{\rm Sl}\nolimits}
\def\SO{\mathop{\rm SO}\nolimits}
\def\span{\mathop{\rm span}\nolimits}

\def\Sp{\mathop{\rm Sp}\nolimits}
\def\Spec{\mathop{\rm Spec}\nolimits}

\def\SS{\mathop{\rm S}\nolimits}

\def\SU{\mathop{\rm SU}\nolimits}
\def\sup{\mathop{\rm sup}\nolimits}

\def\tr{\mathop{\rm tr}\nolimits}
\def\UU{\mathop{\rm U}\nolimits}

\def\0{{\bf 0}}
\def\1{{\bf 1}}

\def\a{{\frak a}}

\def\e{{\frak e}}

\def\g{{\frak g}}

\def\h{{\frak h}}

\def\k{{\frak k}}
\def\l{{\frak l}}

\def\n{{\frak n}}

\def\p{{\frak p}}

\def\r{{\frak r}}
\def\s{{\frak s}}
\def\hsp{{\frak {hsp}}}

\def\sp{{\frak {sp}}}

\def\su{{\frak {su}}}
\def\so{{\frak {so}}}
\def\sL{{\frak {sl}}}
\def\t{{\frak t}}
\def\uu{{\frak u}}

\def\z{{\frak z}}

\def\L{\mathop{\bf L}\nolimits}

\def\C{{\Bbb C}}

\def\N{{\Bbb N}}

\def\R{{\Bbb R}} 
 
\def\Z{{\Bbb Z}} 

\def\:{\colon}  
\def\.{{\cdot}}
\def\|{\Vert}
\def\bsk{\bigskip}

\def\giantskip{\vskip2\bigskipamount}
\def\gsk{\giantskip}
\def \la {\langle}
\def\msk{\medskip}
\def \ra {\rangle}
\def \res {\!\mid\!\!}

\def\ssk{\smallskip}

\def\bbr{\bigbreak}
\def\giantbreak{\par \ifdim\lastskip<2\bigskipamount \removelastskip
         \penalty-400 \giantskip\fi}

\def\nin{\noindent}
\def\cen{\centerline}
\def\pagebreak{\vskip 0pt plus 0.0001fil\break}
\def\linebreak{\break}

\def\hat{\widehat}

\def\epsilon{\varepsilon}
\def\eset{\emptyset}

\def\nin{\noindent}
\def\oline{\overline}

\def\pder#1,#2,#3 { {\partial #1 \over \partial #2}(#3)}
\def\pde#1,#2 { {\partial #1 \over \partial #2}}
\def\phi{\varphi}


\def\subeq{\subseteq}
\def\supeq{\supseteq}

\def\tilde{\widetilde}

\font\eightrm=cmr8


\font\smc=cmcsc10
\font\bfone=cmbx10 scaled\magstep1 
\font\bftwo=cmbx10 scaled\magstep2 

\def\qed{{\unskip\nobreak\hfil\penalty50\hskip .001pt \hbox{}\nobreak\hfil
          \vrule height 1.2ex width 1.1ex depth -.1ex
           \parfillskip=0pt\finalhyphendemerits=0\medbreak}\rm}

\def\qeddis{\eqno{\vrule height 1.2ex width 1.1ex depth -.1ex} $$
                   \medbreak\rm}

\def\Lemma #1. {\bigbreak\vskip-\parskip\noindent{\bf Lemma #1.}\quad\it}

\def\Sublemma #1. {\bigbreak\vskip-\parskip\noindent{\bf Sublemma #1.}\quad\it}

\def\Proposition #1. {\bigbreak\vskip-\parskip\noindent{\bf Proposition #1.}
\quad\it}

\def\Corollary #1. {\bigbreak\vskip-\parskip\nin{\bf Corollary #1.}
\quad\it}

\def\Theorem #1. {\bigbreak\vskip-\parskip\noindent{\bf Theorem #1.}
\quad\it}

\def\Definition #1. {\rm\bigbreak\vskip-\parskip\noindent{\bf Definition #1.}
\quad}

\def\Remark #1. {\rm\bigbreak\vskip-\parskip\noindent{\bf Remark #1.}\quad}

\def\Example #1. {\rm\bigbreak\vskip-\parskip\noindent{\bf Example #1.}\quad}

\def\Problems #1. {\bigbreak\vskip-\parskip\noindent{\bf Problems #1.}\quad}
\def\Problem #1. {\bigbreak\vskip-\parskip\noindent{\bf Problems #1.}\quad}

\def\Conjecture #1. {\bigbreak\vskip-\parskip\noindent{\bf Conjecture #1.}\quad}

\def\Proof#1.{\rm\par\ifdim\lastskip<\bigskipamount\removelastskip\fi\smallskip
            \noindent {\bf Proof.}\quad}

\def\Axiom #1. {\bigbreak\vskip-\parskip\noindent{\bf Axiom #1.}\quad\it}

\def\Satz #1. {\bigbreak\vskip-\parskip\noindent{\bf Satz #1.}\quad\it}

\def\Korollar #1. {\bbr\vskip-\parskip\nin{\bf Korollar #1.} \quad\it}

\def\Bemerkung #1. {\rm\bigbreak\vskip-\parskip\noindent{\bf Bemerkung #1.}
\quad}

\def\Beispiel #1. {\rm\bigbreak\vskip-\parskip\noindent{\bf Beispiel #1.}\quad}
\def\Aufgabe #1. {\rm\bigbreak\vskip-\parskip\noindent{\bf Aufgabe #1.}\quad}

\def\Beweis#1. {\rm\par\ifdim\lastskip<\bigskipamount\removelastskip\fi
           \smallskip\noindent {\bf Beweis.}\quad}

\nopagenumbers

\def\date{\ifcase\month\or January\or February \or March\or April\or May
\or June\or July\or August\or September\or October\or November
\or December\fi\space\number\day, \number\year}

\def\title{Title ??}
\def\author{Author ??}

\def\thanks#1{\footnote*{\eightrm#1}}

\def\rightheadline{\hfil{\eightrm\title}\hfil\tenbf\folio}
\def\leftheadline{\tenbf\folio\hfil{\eightrm\author}\hfil}
\headline={\vbox{\line{\ifodd\pageno\rightheadline\else\leftheadline\fi}}}

\def\firstheadline{}
\def\firstfootline{\cen{\rm\folio}}

\def\seite #1 {\pageno #1
               \headline={\ifnum\pageno=#1 \firstheadline
               \else\ifodd\pageno\rightheadline\else\leftheadline\fi\fi}
               \footline={\ifnum\pageno=#1 \firstfootline\else{}\fi}}

\newdimen\dimenone
 \def\checkleftspace#1#2#3#4{
 \dimenone=\pagetotal
 \advance\dimenone by -\pageshrink   
 \ifdim\dimenone>\pagegoal          
   \else\dimenone=\pagetotal
        \advance\dimenone by \pagestretch
        \ifdim\dimenone<\pagegoal
          \dimenone=\pagetotal
          \advance\dimenone by#1         
          \setbox0=\vbox{#2\parskip=0pt                
                     \hyphenpenalty=10000
                     \rightskip=0pt plus 5em
                     \noindent#3 \vskip#4}    
        \advance\dimenone by\ht0
        \advance\dimenone by 3\baselineskip   
        \ifdim\dimenone>\pagegoal\vfill\eject\fi
          \else\eject\fi\fi}


\def\subheadline #1{\nin\bigbreak\vskip-\lastskip
      \checkleftspace{0.7cm}{\bf}{#1}{\medskipamount}
          \indent\vskip0.7cm\centerline{\bf #1}\medskip}

\def\sectionheadline #1{\bigbreak\vskip-\lastskip
      \checkleftspace{1.1cm}{\bf}{#1}{\bigskipamount}
         \vbox{\vskip1.1cm}\cen{\bfone #1}\bsk}

\def\lsectionheadline #1 #2{\bigbreak\vskip-\lastskip
      \checkleftspace{1.1cm}{\bf}{#1}{\bigskipamount}
         \vbox{\vskip1.1cm}\cen{\bfone #1}\msk \cen{\bfone #2}\bsk}

\def\lchapterheadline #1 #2{\bigbreak\vskip-\lastskip\indent\vskip3cm
                       \cen{\bftwo #1} \msk \cen{\bftwo #2} \gsk}
\def\llsectionheadline #1 #2 #3{\bigbreak\vskip-\lastskip\indent\vskip1.8cm
\cen{\bfone #1} \msk \cen{\bfone #2} \msk \cen{\bfone #3} \nobreak\bsk\nobreak}


\newtoks\literat
\def\[#1 #2\par{\literat={#2\unskip.}%
\hbox{\vtop{\hsize=.15\hsize\nin [#1]\hfill}
\vtop{\hsize=.82\hsize\nin\the\literat}}\par
\vskip.3\baselineskip}

\mathchardef\emptyset="001F 
\def\address{Author: \tt$\backslash$def$\backslash$address$\{$??$\}$}

\def\firstpage{\nin
{\obeylines \parindent 0pt }
\vskip2cm
\centerline {\bfone \title}
\gsk
\centerline{\bf\author}

\vskip1.5cm \rm}

\def\addresstwo{}

\def\dlastpage{\par\vbox{\vskip1cm\nin
\line{
\vtop{\hsize=.5\hsize{\parindent=0pt\baselineskip=10pt\nin\address}}
\quad 
\vtop{\hsize=.42\hsize\nin{\parindent=0pt
\baselineskip=10pt\addresstwo}}
\hfill} }}


\pageno=1

\def\title{Discrete group actions on Stein domains in complex Lie groups}
\def\author{ Dehbia Achab, Frank Betten${}^{\dagger}$, 
Bernhard Kr\"otz${}^*$}
\footnote{}{${}^{\dagger}$Supported by the DFG-grant BE 2080/1--1, 
${}^*$Supported by the DFG-project HI 412/5-2}
\def\date{March  15, 2001}
\def\leftheadline{\tenbf\folio\hfil\eightrm\date}
\def\Box #1 { \msk\par\nin 
\centerline{
\vbox{\offinterlineskip
\hrule
\hbox{\vrule\strut\hskip1ex\hfil{\smc#1}\hfill\hskip1ex}
\hrule}\vrule}\msk }
\def\bs{\backslash}

\def\address
{Dehbia Achab 

Universit\'e Pierre et Marie Curie

Institut de Math\'ematiques 

Analyse Alg\'ebrique - Case 82

Tour 46-0

4, place Jussieu 

F--75252 Paris Cedex 05

France

\gsk\nin

Bernhard Kr\"otz

The Ohio State University 

Department of Mathematics

231 West 18th Avenue 

Columbus, OH 43210-1174

USA

}

\def\addresstwo      
{Frank Betten

Mathematisches Institut

Universit\"at G\"ottingen

Bunsenstra\3e 3--5

D--37073 G\"ottingen

Germany}

\firstpage 

\subheadline{Abstract}

This paper deals with the analytic continuation of holomorphic 
automorphic forms on a (Hermitian) Lie group $G$. We prove that for any 
discrete subgroup $\Gamma$ of $G$ there always exists a 
non-trivial holomorphic automorphic form, i.e., there exists a $\Gamma$-spherical 
unitary highest weight representation of $G$. 
Holomorphic automorphic forms have the property 
that they analytically extend  to holomorphic functions on a complex 
Ol'shanski\u\i{} semigroup $S\subeq G_\C$. 
As an application we prove that the bounded
holomorphic functions on 
$\Gamma\bs S\subseteq \Gamma\bs G_\C$ separate the points.

\sectionheadline{Introduction} 
 Let $G$ be a connected Lie group sitting in its complexification $G_\C$. Let $\Gamma< G$ be 
any discrete subgroup. Then $\Gamma$ acts freely and properly
discontinuously on $G_\C$ and the quotient $\Gamma\bs G_\C$ is a complex manifold. 
It is reasonable to ask which such homogeneous complex manifolds are Stein. 
If $G$ is nilpotent, then $\Gamma\bs G_\C$ is 
Stein by a theorem of Gilligan and Huckleberry (cf.\ [GiHu78]). 
We call a Lie algebra $\g$ {\it weakly elliptic},  if all operators 
$\ad (X)$, $X\in\g$, have imaginary spectrum and note that all nilpotent and many 
solvable Lie algebras are weakly elliptic. Then Loeb has shown in 
[Lo84] that $\Gamma\bs G_\C$ is Stein whenever the Lie algebra $\g$ of $G$ is weakly elliptic. 
But if $G$ is semisimple and $\Gamma$ is a lattice, then a theorem of Barth and Otte 
(cf.\  [BaOt73]) implies that the only holomorphic functions on $\Gamma\bs G_\C$ are the 
constants; in particular $\Gamma\bs G_\C$ is not Stein. 
Therefore in general the $G\times G$-biinvariant complex 
manifold $G_\C$ is too large for $\Gamma\bs G_\C$ carrying a rich complex structure. 
In this paper we will identify an open $G\times G$-invariant Stein  domain 
$S\subeq G_\C$ for which we prove that there is a rich supply of 
automorphic forms on $\Gamma\bs G$ which 
holomorphically extend to $\Gamma\bs S$. In an important example we 
will prove that $\Gamma\bs S$ is Stein.

\par Let us illustrate the situation for the example $G=\SU(p,q)$. Then $G$ can be defined as 
the invariance group of a hermitian form $\la\cdot, \cdot\ra_{p,q}$ on $\C^n$, $n=p+q$, 
and we have $G_\C=\Sl(n,\C)$. The contraction semigroup of the form 
$\la\cdot, \cdot\ra_{p,q}$  
$$S\:=\{g\in G_\C\: \la g.v, g.v\ra_{p,q} <\la v, v\ra_{p,q} \ 
\hbox{for all}\  v\in \C^n, \ v\neq 0\}$$
is easily seen to be an open $G\times G$-invariant domain
in $G_\C$, a so-called {\it complex Ol'shanski\u\i{} semigroup}. 
In general one can define complex  Ol'shanski\u\i{} semigroups for every connected Lie group,  
whenever its Lie algebra admits an open convex weakly elliptic 
$\Ad(G)$-invariant cone $\eset\neq W\subeq \g$. If $\g$ is simple, this means that 
$\g$ is either hermitian or compact, but we also have  many non-reductive examples
as the Jacobi-algebra $\h_n\rtimes\sp(n,\R)$ with $\h_n$ the $2n+1$-dimensional 
Heisenberg algebra (cf.\ [Ne99a] and Example \ 1.3 below). The complex Ol'shanski\u\i{} semigroup 
associated to $G$ and $W$ is given by $S=G\Exp(iW)$. 

\par  In [Ne98, Ne99b] Neeb settled a conjecture of Gindikin by showing 
that all complex 
Ol'shanski\u\i{} semigroups and their symmetric space analogues are Stein 
manifolds. Let now $\Gamma<G$ be any discrete subgroup. In this paper 
we investigate the quotients $\Gamma\bs S$. These manifolds were first considered 
by the first author for $G=\Sl(2,\R)$ in the context of Hardy spaces on
$\Gamma\bs S$ (cf.\ [Ach99]).

Let us in the following assume that $\Ad(G)$ is closed in $\Aut(\g)$ and that 
$W$ is pointed, which both are very natural assumptions in our context. 
Then one of our main results is:

\msk 

\nin {\bf Theorem A.} {\it (cf.\ 4.7)   If $\Gamma<G$ is a discrete subgroup of $G$, 
then: 
\item{\rm (i)} $\Hol(\Gamma\bs S)$ separates the points of $\Gamma\bs S$. 
\item{\rm (ii)} If $\Gamma<G$ is cocompact, then the bounded holomorphic 
functions on $\Gamma\bs S$ separate 
the points. In particular, the Caratheodory semimetric on $\Gamma\bs S$ is 
a metric.}
\msk 

\par We construct holomorphic functions on $\Gamma\bs S$ with techniques from representation theory.  

\par Let $\t\subeq \g$ denote a compactly embedded Cartan subalgebra of $\g$. 
If $(\pi_\lambda, {\cal H}_\lambda)$ is a unitary highest weight representation of $G$
with highest weight $\lambda\in i\t^*$, then one has 
a holomorphic extension to a representation of the complex Ol'shanski\u\i{} semigroup 
$S$ which has $G$ as some sort of Shilov boundary (cf.\ [Ne99a, Ch.\ XIV]). 
We obtain holomorphic 
functions on $\Gamma\bs S$ by taking matrix-coefficients with $\Gamma$-fixed 
hyperfunction vectors; i.e., we take 
$$\theta_{v,\eta}\: \Gamma\bs S\to\C,\ \  \Gamma s\mapsto \la \pi_\lambda(s).v, \eta\ra \, ,$$
with $v\in {\cal H}_\lambda$ and $\eta\in {\cal H}_\lambda^{-\omega}$ a $\Gamma$-fixed element. Here 
${\cal H}_\lambda^{-\omega}$ denotes the Fr\'echet $G$-module of hyperfunction vectors of $(\pi_\lambda, 
{\cal H}_\lambda)$, 
the strong antidual of the analytic vectors ${\cal H}_\lambda^\omega$.   
Note that for $v$ being a highest weight vector the function $\theta_{v,\eta}$ is an analytic continuation of 
the (anti)holomorphic automorphic form $\theta_{v,\eta}\res_G$ (cf.\ [Bo66]).  

\par Our considerations from above show that one gets many holomorphic functions on $\Gamma\bs S$ provided 
$\Gamma$-spherical unitary highest weight representations $(\pi_\lambda, {\cal H}_\lambda)$ 
of $G$ do exist. Our central result is as follows: 

\msk
\nin {\bf Theorem B.} {\it (cf.\  3.8 - 3.11) Let $\Gamma$ be an arbitrary discrete subgroup of $G$
and assume that $G$ has compact center. Denote by $\Lambda\subeq i\t^*$ the spectrum of the 
$L^1$-Bergman space ${\cal B}^1(S)$ on $S$. Then the following assertions hold: 
\item{(i)} For all $\lambda\in \Lambda$ and $v\in {\cal H}_\lambda^\omega$ the Poincar\'e series 
of $v$ 
$$P(v)\:=\sum_{\gamma\in \Gamma} \pi_\lambda(\gamma).v$$
converges in the module of hyperfunction vectors ${\cal H}_\lambda^{-\omega}$ to a 
$\Gamma$-fixed element. 
\item{(ii)} For all but finitely many $\lambda\in \Lambda$ we have $P(v_\lambda)\neq 0$, where 
$v_\lambda$ is a highest weight vector of $(\pi_\lambda, {\cal H}_\lambda)$.}

\msk 
Part (i) is obtained from integral geometric observations in the Bergman space ${\cal B}^1(S)$
together with the Plancherel theorem for ${\cal B}^2(S)$ (cf.\ [Kr98a]). 
\par To get a feeling for the contents of (ii) we explain it for the example $G=\SU(1,1)$. 
Write $D\cong G/K$ for the open unit disc and write elements of $\Gamma$ as 
$\gamma=\pmatrix{a_\gamma& b_\gamma\cr c_\gamma & d_\gamma\cr}$. Then ${\cal H}_\lambda^{-\omega}$
can be naturally be identified with $\Hol (D)$ and we have 
$$P(v_\lambda)(z)\:=\sum_{\gamma\in \Gamma} (c_\gamma z+ d_\gamma)^\lambda \qquad (z\in D).\leqno(*)$$
Here the parameter $-\lambda$ ranges in the  integers with $\lambda\geq 3$. The idea 
for the proof lies in an analytic continuation of $P(v_\lambda)$ in the parameter 
$\lambda$ together with a clever argument that this analytic continuation is non-zero.  
\par Statements related to (ii) were also proved by Siegel (cf. [Fr83]). 
He considers scalar valued Poincar\'e series (this means that the minimal $K$-type of $(\pi_\lambda, {\cal H}_\lambda)$
is one-dimensional) and proves that for a given $\Gamma$ there always exists a $K$-finite vector 
$v$ such that $P(v)\neq 0$. Siegel's proof relies heavily on the fact that $D$ is a bounded domain
while in our approach this is immaterial (for the Jacobi group $D$ is unbounded). Further 
our argument does not need that the Poincar\'e series is scalar valued. 
Finally, Siegel's method does not tell for which $K$-finite vector $v$ one has 
$P(v)\neq 0$ while our techniques yield $P(v_\lambda)\neq 0$ for the distinguished 
highest weight vectors $v_\lambda$. 
We would also like to mention that 
(ii) partly solves an open problem posed by Wells and Wolf (cf.\ [WeWo79]; see also [WaWo83]). 

\msk In Chapter  4 we show that the analytically 
continued automorphic functions $\theta_{v,\eta}$ vanish at infinity: 

\msk \nin {\bf Theorem  C.} {\it (cf.\ 4.6)  If $\Gamma<G$ is cocompact and $v\in {\cal H}_\lambda^\omega$,
then we have 
$$\lim_{s\to \infty \atop s\in \Gamma\bs \oline S} \theta_{v,\eta}(s)=0.$$}

\msk It is our impression that this result in combination with Theorem A 
might yield the Stein property of $\Gamma\bs S$ in the cocompact case (cf.\ Conjecture 5.3 below 
for a more detailed discussion). 
\par For $G=\Sl(2,\R)$ the existence of a non-vanishing holomorphic cusp form (the discriminant)
allows to prove the following: 

\msk\nin {\bf Theorem D.} {\it (cf. \ 5.1 - 5.2)  For $G=\Sl(2,\R)$ and all subgroups $\Gamma<\Sl(2,\Z)$
the quotients $\Gamma\bs S$ are  Stein.}

\msk The authors would like to thank Karl-Hermann Neeb for his very useful remarks.

\sectionheadline{1. Complex Ol'shanski\u\i{} semigroups and their 
representations}

In this first section we introduce complex Ol'shanski\u\i{}
semigroups, which may be thought of as complex Lie
subsemigroups of complex Lie groups. Then we recall the basic facts
concerning their holomorphic representations and the
characterization of the irreducible ones by highest weight
representations.  The concepts will be illustrated with some importnat examples 
(cf.\ Example 1.3 and Remark 1.7).

\subheadline{Complex Ol'shanski\u\i{} semigroups}
Throughout $\g$ denotes a finite dimensional real Lie algebra. 
The concept of a complex Ol'shanski\u\i{} semigroup is closely related
to the concept of elliptic elements in a real Lie algebra $\g$.
Recall that an element $X\in \g$ is called {\it elliptic} if 
$\ad X$ is semsimple with purely imaginary spectrum. Accordingly we 
call a subset $W\subeq \g$ 
{\it elliptic} if all its elements are elliptic.

\Definition 1.1. (Complex Ol'shanski\u\i{} semigroups, cf.\ [Ne99a, Ch.\ XI])
Let $\eset\neq W\subeq \g$ be an open convex elliptic cone and
$\oline W$ its closure.
Let $\tilde G$, respectively $\tilde G_\C$, be the simply connected 
Lie groups associated to $\g$, respectively $\g_\C$, and set 
$G_1\:=\la \exp\g\ra\subeq 
\tilde G_\C$. Then Lawson's Theorem  says that the 
subset $\Gamma_{G_1}(\oline W)\:=G_1\exp(i\oline W)$ is a closed subsemigroup of $\tilde G_\C$ 
and the polar map 
$$G_1\times \oline W\to \Gamma_{G_1}(\oline W), \ \ 
(g, X)\mapsto g\exp(iX)\, ,$$
is a homeomorphism. 
\par Now the universal covering semigroup $\Gamma_{\tilde G}(\oline W)\:=
\tilde\Gamma_{G_1}(\oline W)$ has a similar structure. We can lift the exponential function 
$\exp\: \g + i \oline W\to \Gamma_{G_1}(\oline W)$ to an exponential mapping 
$\Exp\: \g + i\oline W\to \Gamma_{\tilde G}(\oline W)$ with $\Exp(0)=\1$ and thus obtain 
a polar map $\tilde G\times \oline W\to \Gamma_{\tilde G} (\oline W),\ (g, X)\mapsto g\Exp(iX)$
which is a homeomorphism. 
\par If $G$ is any connected Lie group associated to $\g$, then $\pi_1(G)$ is a 
discrete central subgroup of $\Gamma_{\tilde G}(\oline W)$ and we obtain a covering 
homomorphism $\Gamma_{\tilde G}(\oline W)\to \Gamma_G(\oline
W)\:=\Gamma_{\tilde G}(
\oline W)/ \pi_1(G)$. It is easy to see that there is also a polar map 
$G\times \oline W \to \Gamma_G(\oline W), (g, X)\mapsto g\Exp(iX)$, which is a homeomorphism. 
The semigroups of the type $\Gamma_G(\oline W)$ are called {\it complex Ol'shanski\u\i{}
semigroups}.
\par The subset $\Gamma_G(W)\subeq \Gamma_G(\oline W)$ is an open semigroup carrying 
a complex manifold structure such that the multiplication is holomorphic. 
Moreover there is an involution on $\Gamma_G(\oline W)$ given by 
$$^*\: \Gamma_G(\oline W)\to \Gamma_G(\oline W),\ \ s=g\Exp(iX)\mapsto s^*=\Exp(iX)g^{-1} \, ,$$
being antiholomorphic on $\Gamma_G(W)$.
Thus both $\Gamma_G( W)$ and $\Gamma_G(\oline W)$  are involutive 
semigroups.\qed

{}From now on we denote by $S$ an open complex Ol'shanski\u\i{}
semigroup $\Gamma_G(W)$ and by $\oline S$ its ``closure'' $\Gamma_G(\oline
W)$. 

\ssk  If $\l$ is a subalgebra of a Lie algebra $\g$, then we write 
$\Inn_\g(\l)\subeq \Aut(\g)$ for the subalgebra of the automorphism 
group $\Aut(\g)$ of $\g$ which is generated by the elements $e^{\ad X}$, $X\in \l$. 
If $\l=\g$, then we also write $\Inn(\g)$ instead of $\Inn_\g(\g)$. 

\ssk A subalgebra $\a\subeq \g$ is said to be {\it compactly embedded}, 
if $\Inn_\g(\a)$ is relatively compact in $\Aut(\g)$. Note that 
a subalgebra is compactly embedded if and only if it is elliptic. 

\par Note that if a real Lie algebra admits a non-empty open elliptic convex
cone, then there exists a compactly embedded Cartan subalgebra
$\t\subeq \g$ (cf.\ [Ne99a, Th.\ VII.1.8]). To step further we first need some
terminology concerning Lie algebras with compactly embedded Cartan
subalgebras.

\Definition 1.2. (cf.\ [Ne99a, Ch.\ VII]) Let $\g$ be a finite dimensional Lie algebra over
$\R$ with compactly embedded Cartan subalgebra $\t$.
\par\nin (a) Associated to the Cartan subalgebra $\t_\C$ in 
the complexification 
$\g_\C$ there is a root decomposition as follows. For a linear functional 
$\alpha\in \t_\C^*$ we set 
$$\g_\C^\alpha\:=\{X\in \g_\C\: (\forall Y\in \t_\C)\ [Y,X]=\alpha(Y)X\}$$  
and write $\Delta\:=\{\alpha\in\t_\C^*\bs\{0\}\:\g_\C^\alpha\neq \{0\}\}$ for the 
set of roots. Then $\g_\C=\t_\C\oplus\bigoplus_{\alpha\in \Delta}\g_\C^\alpha$, 
$\alpha(\t)\subeq i\R$ for all $\alpha\in \Delta$, and $\oline{\g_\C^\alpha}=
\g_\C^{-\alpha}$, where $X \mapsto \oline X$ denotes complex conjugation on $\g_\C$
with respect to $\g$. 
\par\nin (b)  Let $\r$ denote the radical of $\g$ and note that there is
a $\t$-invariant Levi decomposition $\g=\r\rtimes \s$. Let $\k$ be a maximal compactly embedded
subalgebra of $\g$ containing $\t$. Then a root $\alpha$ is said to be
{\it compact} 
if $\g_\C^\alpha\subeq \k_\C$ and {\it non-compact} otherwise. We write $\Delta_k$ for the set of 
compact roots and $\Delta_n$ for the non-compact ones. 
\par\nin (c) A positive system $\Delta^+$ of roots is a subset of $\Delta$ for which there 
exists  a regular element $X_0\in i\t^*$ with $\Delta^+\:=\{ \alpha\in 
\Delta\: \alpha(X_0)>0\}$. We call a positive system $\k${\it-adapted} if the set 
$\Delta_n^+\:=\Delta_n\cap\Delta^+$ is invariant under the {\it Weyl group}
${\cal W}_\k\:=N_{\Inn_\g(\k)}(\t)/ Z_{\Inn_\g(\k)}(\t)$ acting on $\t$. Recall 
that there exists a $\k$-adapted positive system if and only if 
$\z_\g(\z(\k))=\k$. In this case we say $\g$ is {\it quasihermitian}. In this case it is 
easy to see that $\s$ is quasihermitian, too, and so all simple ideals of $\s$ are either 
compact or hermitian.  
\par\nin (d) We associate to $\Delta_n^+$ the convex cones 
$$C_{\rm min}\:=\oline{\cone\{ i[\oline {X_\alpha},X_\alpha]\: X_\alpha\in \g_\C^\alpha, 
\alpha\in \Delta_n^+\}},$$
and $C_{\rm max}\:=(i\Delta_n^+)^\star=\{X\in\t\: (\forall\alpha\in \Delta_n^+)
\ i\alpha(X)\geq 0\}$. Note that both $C_{\rm min}$ and $C_{\rm max}$ 
are closed 
convex cones in $\t$. 
\par\nin (e) Write $p_\t\:\g\to \t$ for the orthogonal projection along $[\t,\g]$ 
and set ${\cal O}_X\:=\Inn(\g).X$ for the adjoint orbit through $X\in \g$.
We define the {\it minimal} and {\it maximal cone} associated to $\Delta_n^+$ by 
$$W_{\rm min}\:=\{X\in \g\: p_\t({\cal O}_X)\subeq C_{\rm min}\}\quad\hbox{and}
\quad W_{\rm max}\:=\{ X\in \g\: p_\t({\cal O}_X)\subeq C_{\rm max}\}$$   
and note that both cones are convex, closed and $\Inn(\g)$-invariant. \qed

\Example 1.3. (a) (Hermitian case) Let $\g$ be a semisimple Lie algebra
with Cartan decomposition $\g=\k\oplus \p$. Then $\g$ is called {\it hermitian} 
if $\g$ is simple and $\z(\k)\neq \{0\}$. Hermitian Lie algebras are 
classified; the complete list is as follows (cf.\ [Hel78, p.\ 518]):

$$\su(p,q)\quad  \sp(n,\R)\quad  \so^*(2n) \quad \so(2,n)\quad \e_{6(-14)} \quad \e_{7(-25)}.$$
That $\g$ is hermitian implies in particular that there is a compactly embedded 
Cartan subalgebra $\t\subeq \k$ of $\g$  and that $\z(\k)=\R X_0$ is one dimensional (cf.\ [Hel78, Ch.\ VIII]). 
Since $X_0$ is $\Inn_\g(\k)$-fixed, the prescription 
$$\Delta_n^+\:=\{ \alpha\in \Delta\: \alpha(iX_0)>0\}$$
defines a $\k$-adapted system of positive non-compact roots. 
Hermitian Lie algebras admit pointed convex $\Inn(\g)$-invariant closed cones with non-empty 
elliptic interior. One of them, the minimal cone, has a quite simple description: 
$$W_{\rm min}=\oline{\Inn(\g).\R^+X_0}.$$
Its trace with the Cartan algebra $\t$ is given by 

$$C_{\rm min}=\oline {\cone\{-i\check\alpha\: \alpha\in \Delta_n^+\}}$$
where $\check\alpha\in i\t $ denotes the coroot of $\alpha$ (cf.\ [HiNe93, Ch.\ 7]
for all that). 

\par\nin (b) (Compact case) (cf.\ [HiNe93, Ch.\ 7]) Let $\k$ be a compact semisimple Lie algebra 
and set $\g\:=\k\oplus\R$ (direct Lie algebra sum). If $\la\cdot, \cdot\ra$ denotes an
$\Inn(\k)$-innner product on $\k$, then the prescription 

$$W\:=\{ (X,t)\in \g\: \la X,X\ra<t^2, t>0\}$$ 
defines an open elliptic $\Inn(\g)$-invariant convex cone $\g$. The most prominent 
example herof is $\g=\uu(n)$ with 
$$W\:=\{ X\in \uu(n)\: iX \ \hbox{positive definite}\}.$$
On the group level with $G=\UU(n)$ and $G_\C=\Gl(n,\C)$ the Ol'shanski\u\i{} semigroup 
associated to $-W$ is the unit ball in $G_\C$: 

$$S=G\exp(-iW)=\{ g\in \Gl(n,\C)\: \|g\|<1\}$$
with $\|g\|$ the spectral norm.

\par\nin (c) (Mixed case -- The Jacobi algebra) There are interesting examples 
of Lie algebras admitting invariant elliptic convex cones which are neither semisimple 
nor solvable (see Ex.\ 4.9(b) for a non-nilpotent solvable example). 
One of them is the Jacobi algebra 
$$\g\:=\hsp(n,\R)\:=\h_n\rtimes \sp(n,\R)$$
with $\h_n=\R\oplus \R^{2n}$ the $(2n+1)$-dimensional Heisenberg algebra. Then $\z(\g)=\z(\h_n)\cong \R$ 
is one dimensional and a compactly embedded Cartan subalgebra is given by 
$\t=\z(\g)\oplus \t_s$ with $\t_s$ a compactly embedded Cartan subalgebra of $\sp(n,\R)$. 
For the Jacobi algebra $\g$ the maximal cone has an easy description, namely 

$$W_{\rm max}=\R^+\oplus \R^{2n}\oplus W_{\rm max, s}$$
with $W_{\rm max,s}$ a maximal cone in $\sp(n,\R)$. For more details on the Jacobi
algebra we refer to [KN\'O01, Sect.\ II] and [Ne99a, p.\ 400ff, p.\ 700]. For another mixed example 
see also Example 4.9(c) below. \qed

\subheadline{Holomorphic representations}

\Definition 1.4. Let ${\cal H}$ be a Hilbert space and $B({\cal
H})$ the space of bounded operators on it. 
A {\it holomorphic representation} $(\pi,{\cal H})$ 
of a complex Ol'shanski\u\i{} semigroup $S$ is a 
non-degenerate holomorphic semigroup representation $\pi\: S\to B({\cal H})$ such that 
$\pi$ is a vector-valued holomorphic map satisfying 
$\pi(s^*)=\pi(s)^*$ for all $s\in S$.\qed

Now we take a closer look at the irreducible holomorphic
representations  of $S$. They are obtained by
analytic continuation of unitary highest weight representations of
$G$.

{}From now on we assume that $\g$ contains a compactly embedded Cartan subalgebra 
$\t\subeq \g$ and that there exists a non-empty  open elliptic  convex
cone $W\subeq \g$. Then in view of [Ne99a, Th.\ VII.3.8],  there 
exists a $\k$-adapted positive system $\Delta^+$ such that  
$$W_{\rm min}\subeq \oline W\subeq W_{\rm max}$$
holds, $W_{\rm max}^0$ is elliptic, $W_{\rm min}\cap \t= C_{\rm min}$
and $W_{\rm max }\cap\t= C_{\rm max}$. Recall that every 
elliptic $\Ad(G)$-invariant cone $W\subeq \g$ can be 
reconstructed by its intersection with $\t$, i.e., $W=\Ad(G).C$ with $C\:=W\cap\t$.

\Definition 1.5. (Highest weight modules, cf.\ [Ne99a, Ch.\ X])
Let $\Delta^+$ be a positive system. 

\par\nin (a) For a $\g_{\C}$-module $V$ and $\beta \in \t_{\C}^*$
we write $V^\beta := \{ v \in V : (\forall X \in \t_{\C}) 
X.v = \beta(X)v \}$ for the 
{\it weight space of weight} $\beta$ and 
${\cal P}_V = \{ \beta \: V^\beta \not= \{0\} \}$ 
for the set of weights of $V$. 

\par\nin (b) Let $V$ be a $\g_{\C}$-module and $v \in V^\lambda$, $v\neq 0$, a
$\t_\C$-weight vector. 
We say that $v$ is a {\it primitive element of V} (with respect to 
$\Delta^+$) 
if $\g_\C^\alpha.v = \{0\}$ holds for all $\alpha \in \Delta^+$. 

\par\nin (c) A $\g_{\C}$-module $V$ is called a {\it highest weight 
module} with highest weight $\lambda$ (with respect to $\Delta^+$) 
if it is generated by a primitive element of weight $\lambda$. 

\par\nin (d) Let $G$ be a connected Lie group with Lie algebra $\g$. 
We write $T$ for the analytic subgroup of $G$ corresponding to $\t$. 
Let $(\pi, {\cal H})$ be a unitary representation of $G$. A vector 
$v\in{\cal H}$ is called $T${\it -finite} if it is contained in a
finite dimensional $T$-invariant subspace. 
We write ${\cal H}^{T,\omega}$ for the space of analytic 
$T$-finite vectors. 

\par\nin (e) An irreducible unitary representation $(\pi,{\cal H})$ of 
$G$ is called a {\it highest weight representation} with respect to $\Delta^+$ 
and highest weight $\lambda\in i\t^*$ if 
${\cal H}^{T,\omega}$ is 
a highest weight module for $\g_\C$ with respect to $\Delta^+$ and 
highest weight $\lambda$. 
We write $HW(G, \Delta^+)\subset i\t^*$ for the set of highest weights 
corresponding to unitary highest weight representations of $G$ with
respect to $\Delta^+$.\qed

The interplay between irreducible holomorphic representations of $S$ and 
unitary highest weight representations of $G$ is described in the 
following lemma.  
Recall that a cone $W$ in a real vector space $V$ is called {\it pointed} if $W$ contains no affine 
lines. 

\par The central ideas in the following theorem go back to Ol'shanski\u\i{}
(cf.\ [Ol82]) Stanton (cf.\ [St86]).

\Theorem 1.6. Let $S=\Gamma_G(W)$ be a complex Ol'shanski\u\i{}
semigroup and $\Delta^+$ be a $\k$-adapted positive system with 
$C_{\rm min}\subeq \oline C\subeq C_{\rm max}$. Suppose that $W$ is pointed. 
\item{(1)} If $(\pi, {\cal H})$ is an irreducible holomorphic
representation of $S$, then $(\pi, {\cal H})$ extends to a strongly
continuous representation of 
$\oline S$, also denoted by $(\pi, {\cal H})$, such that $\pi\res_G$
is a uniquely determined unitary highest weight representation of $G$ 
with respect to $\Delta^+$. 
\item{(2)} Conversely, if $(\pi_\lambda, {\cal H}_\lambda)$ is a unitary highest
weight representation of $G$ with respect to  $\Delta^+$, then 
$(\pi_\lambda, {\cal H}_\lambda)$
extends to a uniquely determined strongly continuous representation of
$\oline S$, whose restriction to $S$ is holomorphic and irreducible. 

\Proof. This follows from [Ne99a, Th.\ XI.4.8] together with its 
following remark.\qed 

We assume in the following that the cone $W\subeq \g$ is pointed. 

\Remark 1.7. In the special case of the Jacobi group 
$$G\:=\HSp(n,\R)=H_n\rtimes \Sp(n,\R)$$
(a group associated to the Jacobi algebra $\hsp(n,\R)$; cf.\ Example 1.3(c))
it is interesting to see what the unitary highest weight representations of 
$G$ are. They are given by 
$$\sigma_\mu\otimes \pi_{\lambda,s}$$
where $\sigma_\mu$ is an extended metaplectic representation of $G$ and 
$\pi_{\lambda,s}$ is a unitary highest weight representation of $\Sp(n,\R)$
extended trivially to the Heisenberg group $H_n$ (cf.\ [Ne99a, Sect.\ X.3]). 
All these representations extend to holomorphic representations of the maximal
open complex Ol'shanki\u\i{} semigroup $\Gamma_G(\Int W_{\rm max})$.\qed

\sectionheadline{2. The manifolds $\Gamma\bs S $ and some concepts from complex analysis}

Let $S=\Gamma_G(W)=G\Exp(iW)$ be an open complex Ol'shanski\u\i{} semigroup and 
$\Gamma<G$ a discrete subgroup. In view of the polar decomposition of $S$ 
(cf.\ Definition 1.1), 
the group $\Gamma$ acts via 
$$\Gamma\times S\to S, \ \ (\gamma, s)\mapsto \gamma s\, ,$$
freely and properly discontinuously on $S$. Thus the quotient $\Gamma\bs S$ is 
Hausdorff and carries a complex structure which is induced 
from the quotient map $S\to \Gamma\bs S, \ s\mapsto \Gamma s$.

\subheadline{Stein manifolds}

Since one of the main objectives of this paper is the investigation 
of the structure of the complex manifolds $\Gamma\bs S$, especially 
whether they are Stein,  we 
briefly recall here some facts concerning Stein manifolds.

\par If $M$ is a complex manifold, then we write $\Hol(M)$ for the space of
holomorphic functions on $M$. We endow $\Hol(M)$ with the topology of
compact convergence on $M$. Note that $\Hol(M)$ is a Fr\'echet space
provided $M$ is second countable.

\Definition 2.1. (Stein manifolds) Let $M$ be a second countable
complex manifold. Then $M$ is called a {\it Stein manifold}, if the
following two axioms are satisfied

\item {(St1)} The points of $M$ are separated by $\Hol(M)$.
\item{(St2)} The manifold $M$ is {\it holomorphically convex}, i.e.,
for each compact subset $K\subeq M$ the {\it holomorphically convex
hull}
$$\hat K\:=\{ z\in M\: (\forall f\in \Hol(M)) |f(z)|\leq \sup_{w\in
K} |f(w)|\}$$ is compact.\qed

\Theorem 2.2. {\rm (Grauert)} For a second countable complex manifold $M$, the
following assertions are equivalent:

\item{(1)} The manifold $M$ is Stein. 
\item{(2)} There exists a strictly plurisubharmonic exhaustion function
$\phi\: M\to\R$, i.e., $\phi$ is  strictly plurisubharmonic and for
all $r\in \R$ the set $\phi^{-1}(]-\infty, r])$ is a compact subset
of $M$.  

\Proof. [H\"o73, Th.\ 5.2.10].\qed

\Theorem 2.3. If $E \rightarrow M$ is a holomorphic principal bundle,
such that the fiber and the base~$M$ are Stein, then $E$ is Stein.

\Proof. [MaMo60, Th.\ 4].\qed

\Theorem 2.4. {\rm (Neeb)} Let $S$ be a complex Ol'shanski\u\i{}
semigroup. Then there exists a $G\times G$-invariant strictly 
plurisubharmonic  smooth function $\phi\: S\to\R^+$ with 
$$\lim_{s_n\to s}\phi(s_n)=\infty$$
for all $s\in \oline S\bs S$. 

\Proof. [Ne98, Lemma 5.11].\qed

\subheadline{Invariant Hilbert spaces}

We conclude this section by introducing the concept of an invariant Hilbert space 
of holomorphic functions.

\Definition 2.5. Let $M$ be a complex manifold and $G$ a group acting 
on $M$ via 
$$M\times G\to M, \ \ (m,g)\mapsto m.g\, ,$$
by holomorphic automorphisms. Then a Hilbert space 
${\cal H}\subeq \Hol(M)$ is called a {\it $G$-invariant Hilbert space of holomorphic 
functions on $M$} if the following two axioms are satisfied:
\item{(IH1)} The inclusion map ${\cal H}\into \Hol(M)$ is continuous.
\item{(IH2)} We have a unitary representation 
$$\rho\: G\to U({\cal H}), \ \ (\rho(g).f)(z)\:=f(z.g).\qeddis 

Let $(\rho, {\cal H})$ be an invariant Hilbert space. Then by (IH1) all point evaluations 
${\cal H}\to\C, \ f\mapsto f(z)$ are continuous. Hence there exists for 
every $z\in M$ an element $K_z\in {\cal H}$ such that $\la f, K_z\ra =f(z)$ for all 
$f\in {\cal H}$.  The function 
$$K\: M\times M\to\C, \ \ (z,w)\mapsto K(z,w)\:=\la K_w, K_z\ra\, ,$$
is holomorphic in the first, antiholomorphic in the second variable  and $G$-invariant, i.e., 
we have $K(z.g,w.g)=K(z,w)$ for all $z,w\in M$, $g\in G$. 
We call the function $K$ the {\it reproducing kernel} of $(\rho, {\cal H})$. 
For further information on reproducing kernel Hilbert spaces we refer to 
[Ne99a, Ch.\ I-IV].

\sectionheadline{3. Constructing $\Gamma$-spherical representations}

In this section we give a construction of $\Gamma$-spherical highest weight 
representations for arbitrary discrete subgroups $\Gamma<G$. 
It was an idea of Godement to obtain $\Gamma$-spherical representations by averaging 
matrix coefficients of integrable representations of $G$ (cf.\ [Bo66]). 
In our setup however it is much more convenient to use 
the Bergman space ${\cal B}^1(S)$ of integrable holomorphic functions on $S$ 
instead 
of $L^1(G)$. This is because we can naturally realize  highest weight representations with 
sufficiently large parameter in the Hilbert space ${\cal B}^2(S)$ (cf.\ [Kr98a]).

\subheadline{Hyperfunction vectors for unitary representations}

\Definition 3.1. Let $G$ be a Lie group and ${\cal H}$ a Hilbert space.

\par\nin(a) For a unitary representation $(\pi,{\cal H})$ of $G$ we denote by ${\cal H}^\infty$
and ${\cal H}^\omega$ the space of all smooth, respectively  analytic
vectors of $(\pi ,{\cal H})$.
The corresponding strong antiduals are denoted by ${\cal H}^{-\infty }$ and 
${\cal H}^{-\omega}$ and their elements are called {\it distribution}, 
resp.\ {\it hyperfunction vectors} (see [KN\'O97, Appendix] for the definition of the 
topology of ${\cal H}^\omega$). Note that there is a natural chain of 
continuous inclusions
$${\cal H}^\omega
\into {\cal H}^\infty \into {\cal H}\into 
{\cal H}^{-\infty}\into {\cal H}^{-\omega}.$$
The natural extension of $(\pi,{\cal H})$ to a representation on the space 
of hyperfunction vectors is denoted by $(\pi^{-\omega}, {\cal
H}^{-\omega})$ and given explicitly by 
$$ \la \pi^{-\omega}(g).\nu, v \ra := \la \nu, \pi(g^{-1}).v
\ra. $$
Note that $(\pi^{-\omega}, {\cal H}^{-\omega})$ is a
continuous representation of $G$ (cf.\ [KN\'O97, App.]). 

\par\nin (b) Let $\Gamma\subeq G$ be a closed subgroup. 
For a unitary representation $(\pi,{\cal H})$ of $G$ we write 
$({\cal H}^{-\omega})^\Gamma$ for the set of 
all those elements $\nu\in {\cal H}^{-\omega}$ 
satisfying $\pi^{-\omega}(\gamma).\nu=\nu$
for all $\gamma\in \Gamma$. The unitary representation $(\pi,{\cal H})$ is called 
$\Gamma$-{\it spherical} if there exist a cyclic vector 
$\nu\in  ({\cal H}^{-\omega})^\Gamma$.\qed

\Remark 3.2. Even though in general the topology on the space of analytic 
vectors is hard to deal with, one has a quite explicit 
picture for holomorphic representations of complex Ol'shanski\u\i{}    
semigroups.  Let $G$ be a connected Lie group  and $(\pi, {\cal H})$ 
a unitary representation of $G$ which has an extension to a holomorphic 
representation of some  complex  Ol'shanski\u\i{} semigroup 
$S=G\Exp(iW)$. Then we have for each $X \in iW$ 
$${\cal H}^\omega=\bigcup_{t>0}\pi(\Exp(tX)).{\cal H}, $$  
and the topology on ${\cal H}^\omega$ is  the finest locally 
convex topology on ${\cal H}^\omega$
making for all $t>0$ the maps ${\cal H}\to {\cal H}^\omega, \ \ v\mapsto \pi(\Exp(tX)).v$ continuous
(cf.\ [KN\'O97, Appendix]). Note that we have $\pi(S).{\cal H}^{-\omega}\subeq {\cal H}^\omega$ and that 
the action of $S$ on 
${\cal H}^{-\omega}$ is given by 
$$ \la \pi^{-\omega}(s).\nu, v \ra := \la \nu, \pi(s^*).v \ra.\qeddis 

\par A Lie group $G$ is called a {\it (CA)-Lie group} ((CA)=closed adjoint) if $\Ad(G)$ is closed in 
$\Aut(\g)$. Note that most Lie groups are (CA)-Lie groups; in particular, 
all reductive and nilpotent Lie groups have the (CA)-property (cf.\ [Ne99a, VII.1.13--VII.1.15]).

\Lemma 3.3. Let $(\pi_\lambda, {\cal H}_\lambda)$ be a unitary highest weight representation 
of the (CA)-Lie group $G$ with highest weight $\lambda\in i\Int C_{\rm min}^\star$. Then there exists
a $X\in iC$ with $\Spec ( d\pi_\lambda(X))\subeq -\beta-\N_0$ for some  $\beta\in \R$.
Let $(v_n)_{n\in\N}$ be an orthonormal basis of ${\cal H}_\lambda$ consisting 
of $\t_\C$-weight vectors. 
Then for all $t>0$ there exists a constant $C_t>0$ such that 
$$(\forall w\in {\cal H_\lambda}) \qquad \sum_{n=1}^\infty |\la \pi_\lambda(\Exp(tX)).w,v_n\ra |
\leq C_t \|w\|.$$

\Proof. Let ${\cal P}_\lambda$ denote the $\t_\C$-weights of ${\cal H}_\lambda$. Then 
${\cal P}_\lambda\subeq \lambda-\N_0[\Delta^+]$ as $(\pi_\lambda, {\cal H}_\lambda)$
is a highest weight module with highest weight $\lambda$. 
\par Since $C\subeq \Int C_{\rm max}$ we have $\alpha(iY)>0$ for all 
$\alpha\in\Delta_n^+$ and $Y\in C$. Moreover $C$ is ${\cal W}_\k$-invariant so that we can find an 
element $Y\in C$ such that $\alpha(iY)>0$ for all $\alpha\in \Delta^+$. Since $G$ is a (CA)-Lie group 
the group $T/ Z(G)$ is compact (cf.\ [Ne99a, Cor.\ VII.1.5]) and hence we can choose $X\in iC$ such that 
$\alpha(iX)\in \N$ for all $\alpha\in \Delta^+$. 
This proves  $\Spec ( d\pi_\lambda(X))\subeq -\beta -\N_0$ for 
some $\beta\:=-\lambda(X)\in\R$. 

\par Let ${\cal H}_\lambda=\bigoplus_{n=0}^\infty {\cal H}_\lambda^{-\beta -n}$ be the 
orthogonal decomposition of ${\cal H}_\lambda$ into $d\pi_\lambda(X)$-eigenspaces. 
Set $m_n\:=\dim {\cal H}_\lambda^{-\beta -n}$ and note that 
there exists an $N\in \N$ such that 
$$(\forall n\in\N)\qquad m_n\leq (n+1)^N$$
(cf.\ [Ne99a, Lemma X.4.9]). 
Then 
$(v_n)_{n\in\N}=(v_j^n)_{1\leq j\leq m_n,  n\in \N}$ 
with 
$v_j^n\in {\cal H}_\lambda^{-\beta -n}$ for all $1\leq j\leq m_n$. 

\par Take now $v=\sum_n c_n v_n=\sum_{j,n} c_j^n v_j^n =\pi_\lambda(\Exp(tX)).w$
for some $w\in {\cal H}_\lambda$ and $t>0$. Then 
$w=\sum_{j,n} e^{t(n+\beta)} c_j^n v_j^n$ and we have 
$$\eqalign{\sum_{n=1}^\infty |c_n|&=\sum_{n=1}^\infty \sum_{1\leq j\leq m_n} |c_j^n|=
\sum_{n=1}^\infty \sum_{1\leq j\leq m_n} |c_j^n| e^{t(n+\beta)}e^{-t(n+\beta)}\cr 
&\leq \Big(\sum_{n=1}^\infty \sum_{1\leq j\leq m_n} |c_j^n|^2 e^{2t(n+\beta)}\Big)^{1\over 2}
\Big(\sum_{n=1}^\infty \sum_{1\leq j\leq m_n} e^{-2t(n+\beta)}\Big)^{1\over 2}\cr 
&= \|w\|\Big(\sum_{n=1}^\infty m_n e^{-2t(n+\beta)}\Big)^{1\over 2}\cr 
&\leq \|w\|\Big(\sum_{n=1}^\infty (1+n)^N e^{-2t(n+\beta)}\Big)^{1\over 2}\leq 
C_t\|w\|\cr}$$
with $C_t\:=\Big(\sum_{n=1}^\infty (1+n)^N e^{-2t(n+\beta)}\Big)^{1\over 2}<\infty.$ 
\qed

\subheadline{Bergman spaces}

If $G$ is a locally compact group, then we denote by $\mu_G$  a 
left Haar measure on $G$. 

\par Let $S=\Gamma_G(W)$  be a complex Ol'shanski\u\i{} semigroup and $H<G$ a
unimodular subgroup. Recall from [Kr98a, Sect.\ II] that there exists 
an $\oline S$-invariant positive measure $\mu_{H\bs S}$ on $H\bs S$. 
Denote by $\Hol(S)^H$ the space of left $H$-invariant holomorphic 
functions on $S$. For each $1\leq p\leq \infty$ we write $\|\cdot\|_p$
for the norm of the Banach space $L^p(H\bs S, \mu_{H\bs S})$. For 
$1\leq p\leq \infty$ we now define 
$${\cal B}^p(H\bs S)\:=\{ f\in \Hol(S)^H\: \|f\|_p<\infty\}$$
and note that ${\cal B}^p(H\bs S)$ is a Banach space of holomorphic
functions on $S$ (cf.\ [Kr98a, Prop.\ II.4]); if $p=2$, then 
${\cal B}^2(H\bs S)$ is even a Hilbert space. If $H=\Gamma$ is a
discrete subgroup of $G$, then we have $\Hol(S)^\Gamma\cong \Hol(\Gamma\bs S)$ and 
we call ${\cal B}^2(\Gamma\bs S)$ the 
{\it Bergman space} of the complex manifold $\Gamma\bs S$.

\Proposition 3.4. Let $S$ be a
complex Ol'shanski\u\i{} semigroup and $H<G$ a unimodular 
subgroup. Then the following assertion holds: 
\item{(i)} For all compact subsets $Q\subeq S$ there exists 
a constant $C_Q>0$ such that 
$$(\forall f\in {\cal B}^1(S))(\forall s\in Q)\ 
\int_H |f(hs)|\ d\mu_H(h)\leq C_Q\|f\|_1$$
holds. 
\item {(ii)} The linear map 
$$I^H\: {\cal B}^1(S)\to {\cal B}^1(H\bs S), \ \ f\mapsto f^H; \ 
f^H(s)=\int_H f(hs)\ d\mu_H(h)\, ,$$
is defined and a contraction of Banach spaces. 

\Proof. (i) Let $Q\subeq S$ be a compact subset. W.l.o.g.\ we may
assume that there exists an open relatively compact subset $U\subeq S$
such that $Q\subeq U\subeq \oline U\subeq S$ and that $H\oline U \cong
(H\bs H\oline U)\times H$ holds. Thus we may assume that there exists a compact 
subset $K_H\subeq H$ such that $\oline U\cong (H\bs H\oline U)\times
K_H$. It follows as in [Kr98a, Prop.\ II.4(i)] that there exists a constant $c_Q>0$
such that 
$$(\forall s\in Q)\ |f(s)|\leq c_Q \int_U |f(u)|\
d\mu_S(u)$$
holds for all $f\in \Hol(S)$. For each $f\in \Hol(S)$ and $g\in G$
we set $f_g(s)=f(gs)$. Note that $f_g\in {\cal B}^1(S)$ provided $f\in
{\cal B}^1(S)$ since $\mu_S$ is $G\times G$-biinvariant 
(cf.\ [Kr98a, Sect.\ II]). Thus we get for all $f\in {\cal B}^1(S)$
and $s\in Q$ that 
$$\eqalign{\int_H |f(hs)|\ d\mu_H(h)&=\int_H |f_h(s)|\ d\mu_H(h)\leq
c_Q\int_H \int_U |f_h(u)|\ \ d\mu_S(u)\ d\mu_H(h)\cr
& \leq c_Q\int_H \int_U |f(hu)|\ \ d\mu_S(u)\ d\mu_H(h)\cr
&=c_Q\int_U \int_H |f(hu)|\ \ d\mu_H(h)\ d\mu_S(u)\cr
&=c_Q\int_{H\bs HU}\int_{K_H}  \int_H |f(h_2h_1u)|\ d\mu_H(h_2)\ d\mu_H(h_1)\ d\mu_{H\bs S}(Hu)\cr
&\leq c_Q\mu_H(K_H)\int_{H\bs HU}\int_H |f(hu)|\ d\mu_H(h)\ d\mu_{H\bs S}(Hu)\cr
&\leq c_Q\mu_H(K_H)\int_S|f(s)|\ d\mu_S(s)\leq C_Q\|f\|_1,\cr}$$
where $C_Q=c_Q\mu_H(K_H)$. This proves (i). 

\par\nin (ii) Let $f\in {\cal B}^1(S)$. Then the inequality in (i)
implies that $f^H\in \Hol(S)^H$. It remains to show that
$\|f^H\|_1\leq \|f\|_1$. But this follows from 
$$\eqalign{\|f^H\|_1&=\int_{H\bs S} |f^H(s)|\ d\mu_{H\bs S}=\int_{H\bs S}\Big|\int_H f(hs)\
d\mu_H(h)\Big|\ d\mu_{H\bs S}(Hs)\cr
&\leq \int_{H\bs S}\int_H |f(hs)|\ d\mu_H(h)\ d\mu_{H\bs S}(Hs)=\|f\|_1,\cr}$$
completing the proof of (ii). \qed

Let $(\pi,{\cal H})$ be a unitary representation of a Lie group  $G$ and $(\pi^*, {\cal H}^*)$ its 
dual representation. We denote by  $(\hat\pi, B_2({\cal H}))$ the unitary 
representation of $G\times G$ on Hilbert-Schmidt operators on ${\cal H}$ given by $\hat\pi(g_1,g_2).A=\pi(g_1)A\pi(g_2)^{-1}$
for $g_1,g_2\in G$ and $A\in B_2({\cal H})$. Recall that the scalar product on 
$B_2({\cal H})$ is given by 
$\la A,B\ra =\tr (AB^*)$ for all $A,B\in B_2({\cal H})$. 
Further we have a natural equivalence 
$$(\pi\otimes\pi^*, {\cal H}\hat\otimes {\cal H}^*) \to 
(\hat\pi, B_2({\cal H})), 
\ \ v\otimes \la\cdot,  w\ra \mapsto P_{v,w}; \ P_{v,w}(u)\:=\la u,w\ra v.$$

\par Let $(R\otimes L, \Hol(S))$ denote the right-left regular 
representation of 
$G\times G$ on $\Hol(S)$, i.e., we have 
$$(\forall g_1,g_2\in G)(\forall f\in \Hol(S))(\forall s\in S)\quad 
\big((R\otimes L)(g_1,g_2).f\big)(s)=f(g_2^{-1}sg_1).$$
Recall that for a unitary highest weight representation $(\pi_\lambda, {\cal H}_\lambda)$ all 
operators $\pi_\lambda(s)$, $s\in S$, are of trace class (cf.\ [Ne99a, Th.\ XI.6.1]). 
Further the prescription 
$$(\hat\pi_\lambda, B_2({\cal H}_\lambda))\to (R\otimes L,  \Hol(S)), \ \ A\mapsto 
f_A^\lambda(s)=\tr (A\pi_\lambda(s)).$$
gives us a $G\times G$-equivariant realization  of $B_2({\cal H}_\lambda)$ as 
$G\times G$-invariant Hilbert space of holomorphic functions on $S$ (cf.\ Definition 2.5)

\Theorem 3.5. Assume that $G$ is a (CA)-Lie group and that $Z=Z(G)$ is compact. 
If $S$ is a complex Ol'shanski\u\i {} semigroup, then we have 
a $G\times G$-equivariant isomorphism of Hilbert spaces
$$\hat \bigoplus_{\lambda\in \Lambda}(\hat\pi_\lambda, B_2({\cal H}_\lambda))
\to 
(R\otimes L, {\cal B}^2(S)), \ \ B_2({\cal H}_\lambda)\ni A\mapsto \sqrt{b(\lambda)} f_A^\lambda,$$
with $b(\lambda)>0$ for all $\lambda\in \Lambda$ and 
$$\Lambda\:=\{\lambda\in HW(G, \Delta^+)\cap iC^\star\: \lambda+\rho \in i\Int C_{\rm min}^\star, 
\lambda+2\rho\in i \Int (C\cap\t^+)^\star\}, $$
where $\rho\:={1\over 2}\sum_{\alpha\in \Delta^+}\alpha$ and 
$\t^+\:=\{X\in\t\: (\forall \alpha\in\Delta^+) i\alpha(X)\geq 0\}$. 

\Proof.  This is a special case of [Kr98a, Thm.\ IV.5].\qed 

For the rest of this section we assume now that $G$ is a (CA)-Lie 
group, $Z(G)$ is compact. 
We denote by $K$ the analytic subgroup of $G$ corresponding to $\k$. Note that 
$K$ is compact by our assumptions on $G$.

\Lemma 3.6.  We have $2\Lambda\subeq \Lambda$. 
Further for all $\lambda\in 2\Lambda$ and all $T$-finite analytic vectors $w\in {\cal H}_\lambda^{T,\omega}$
we have a
continuous $G$-equivariant inclusion mapping 
$$\Phi\: {\cal H}_\lambda^\omega \to {\cal B}^1(S), \ \ v\mapsto f_{P_{v,w}}^\lambda\, ,$$
with $f_{P_{v,w}}^\lambda(s)=\la \pi_\lambda(s).v,w\ra$, $s\in S$. 

\Proof. We first show that $2\Lambda\subeq \Lambda$. Recall the description of 
$\Lambda$ in Theorem 3.5. Let $\lambda\in \Lambda$. 
Then $\rho \in -iC_{\rm min}^\star$ implies that 
$$2\lambda+\rho =2(\lambda+\rho)-\rho \in i\Int C_{\rm min}^\star+ 
iC_{\rm min}^\star\subeq  i\Int C_{\rm min}^\star.$$ 
Thus $2\lambda$ is a parameter of the holomorphic discrete series and 
in particular we have $2\lambda\in HW(G,\Delta^+)$ (cf.\ [Ne99a, Sect.\ XII.5] for all that). 
Finally we have 
$$2\lambda+2\rho=(\lambda+2\rho)+\lambda \in i\Int (C\cap\t^+)^\star+ iC^\star\subeq  
i \Int (C\cap\t^+)^\star$$ 
since $\lambda\in i C^\star $ and $C^\star \subeq (C\cap\t^+)^\star$. This 
concludes the proof of $2\Lambda\subeq \Lambda$. 

\par Let now $\lambda\in \Lambda$ and write 
$v_{\lambda}$, respectively $v_{2\lambda}$, for a normalized highest
weight vector in ${\cal H}_\lambda$, respectively
${\cal H}_{2\lambda}$. Denote by 
${\cal H}_\lambda\hat\otimes{\cal H}_\lambda$ the $G$-module where $G$ acts diagonally. Then 
the prescription $v_{2\lambda}\mapsto v_{\lambda} \otimes v_{\lambda}$ gives rise to
an equivariant embedding of ${\cal H}_{2\lambda}$ into ${\cal H}_\lambda\hat\otimes{\cal H}_\lambda$
(for more details see the proof of [HiKr99, Th.\ 3.2.1] where a similar 
situation is dealt with).  
In the sequel we assume that ${\cal H}_{2\lambda}\subeq{\cal H}_\lambda\hat\otimes{\cal H}_\lambda$.
Note that ${\cal H}_{2\lambda}^{T,\omega}\subeq {\cal H}_\lambda\otimes{\cal H}_\lambda$ 
so that by linearity we may 
assume that $w=w_1\otimes w_2$ for some $w_1, w_2\in {\cal H}_\lambda$. 

\par Let $(v_n)_{n\in\N}$ be an orthonormal basis of ${\cal H}_\lambda$. Then
$(v_n\otimes v_m)_{n,m\in\N}$ constitutes an orthonormal basis of ${\cal H}_\lambda\hat\otimes
{\cal H}_\lambda$. Hence $v\in {\cal H}_{2\lambda}^{\omega}$ can be written as 
$v=\sum_{n,m} c_{n,m} v_n\otimes v_m$ with 
$\|v\|^2=\sum_{n,m}|c_{n,m}|^2$.  
Let $X\in iC $ be as in Lemma 3.3. Then $(X,X)\in i(C\times C)$ meets the same assumptions
as $X$ in Lemma 3.3 for the unitary highest weight representation $\pi_\lambda\otimes\pi_\lambda$ of 
$G\times G$. Hence if $v=\pi_{2\lambda}(\Exp(tX)).u$ for some $u\in {\cal H}_{2\lambda}$ then 
$\sum_{n,m} |c_{n,m}|\leq C_t \|u\|$ by Lemma 3.3. 
Now the 
Cauchy-Schwarz inequality together with Theorem 3.5 gives  
$$\eqalign{\int_S |f_{P_{v,w}}^{2\lambda} (s)|&\ d\mu_S(s)\leq \sum_{n,m} |c_{n,m}|
\int_S |\la (\pi_\lambda(s)\otimes\pi_\lambda(s)).(v_n\otimes v_m), w_1\otimes w_2\ra | \ d\mu_S(s)\cr 
&\leq\sum_{n,m} |c_{n,m}| \int_S |\la \pi_\lambda(s).v_n, w_1\ra \la \pi_\lambda(s).v_m, w_2\ra |
\ d\mu_S(s)\cr 
&\leq\sum_{n,m} |c_{n,m}| \Big(\int_S |\la \pi_\lambda(s).v_n, w_1\ra|^2\ d\mu_S(s)\Big)
^{1\over 2}\Big(\int_S |\la \pi_\lambda(s).v_m, w_2\ra|^2\ d\mu_S(s)\Big)
^{1\over 2} \cr 
&={1\over b(\lambda)} \|w_1\|\cdot \|w_2\|\sum_{n,m} |c_{n,m}|\leq c_t \|u\|\cr}$$
with $c_t\:={C_t\over b(\lambda)} \|w_1\|\cdot \|w_2\|<\infty$. 
Thus we see that for all $t>0$ the mapping 
$$ {\cal H}_{2\lambda}^\omega \to {\cal B}^1(S), \ \ 
v \mapsto \Phi(\Exp(tX).v) \, ,$$
is continuous. In view of the definition of the topology of 
${\cal H}_{2\lambda}^{\omega}$, this proves 
the continuity of $\Phi$. \qed

\subheadline{Poincar\'e series}

Now we are going to construct $\Gamma$-spherical highest weight 
representations for generic parameters in $2\Lambda$.

\Definition 3.7. Let $(\pi, {\cal H})$ be a unitary highest weight representation for 
$G$. For $v\in {\cal H}$ we call 
$$P(v)\:=\sum_{\gamma\in \Gamma}\pi(\gamma).v$$
the {\it Poincar\'e series of $v$}, provided $P(v)$ exists
in ${\cal H}^{-\omega}$, i.e., if 
$${\cal H}^\omega\to \C, \ \ w\mapsto \sum_{\gamma\in \Gamma}\la \pi(\gamma).v,w\ra\, ,$$
is a continuous antilinear functional. Note that $P(v)\in ({\cal H}^{-\omega})^\Gamma$.\qed

\Proposition 3.8. Let $(\pi_\lambda, {\cal H}_\lambda)$ be a unitary 
highest weight representation of $G$ with $\lambda\in 2\Lambda$. Then for all 
$w\in {\cal H}_\lambda^{T,\omega}$ the Poincar\'e series $P(w)$ exists,  and defines a 
$\Gamma$-fixed hyperfunction vector. Moreover,  for all $v\in {\cal H}_\lambda^\omega$ the series 
$\la P(w),v\ra = \sum_{\gamma\in \Gamma}\la \pi(\gamma).w,v\ra$ converges 
absolutely.

\Proof. Let $w\in {\cal H}_\lambda^{T,\omega}$,  $v\in {\cal H}_\lambda^\omega$
and set 
$$f_v\: \oline S\to \C, \  \ s\mapsto\la\pi_\lambda(s).v,w\ra.$$
By Lemma 3.6 we have $f_v\in {\cal B}^1(S)$ and so 
$f_v^\Gamma\in {\cal B}^1(\Gamma\bs S)$ with $f_v^\Gamma$ converging absolutely 
on all compact subsets $Q\subeq S$, i.e., 

$$(\exists C_Q>0)(\forall s\in Q)(\forall f\in {\cal B}^1(S))\qquad 
\sum_{\gamma\in \Gamma} |f(\gamma s)|\leq 
C_Q\|f\|_1\leqno (3.1)$$
(cf.\ Proposition 3.4). Let $X\in iC$ and set $s_t\:=\Exp(tX)$ for all $t>0$. Then 
Lemma 3.6 asserts in particular that 
for all $t>0$ there exists a constant $C_t>0$ such that 
$$(\forall v\in {\cal H}_\lambda)\qquad \|f_{\pi_\lambda(s_t).v}\|_1\leq C_t \|v\|.\leqno(3.2)$$   
\par By the definition 
of the topology on ${\cal H}_\lambda^{\omega}$, we only have to check  
the continuity of the maps 
$${\cal H}_\lambda\to\C, \ \ v\mapsto \la \pi_\lambda(s_t).v, P(w)\ra.$$
Fix now $t>0$ and choose $Q$ such that $s\:=s_{t\over 2}\in Q$. 
Then (3.1) and (3.2) imply that 
$$\eqalign{|\la \pi_\lambda(s_t).v, P(w)\ra|&=
|\sum_{\gamma\in \Gamma}\la \pi_\lambda(s_t).v, \pi_\lambda(\gamma).w\ra|\cr 
&\leq\sum_{\gamma\in \Gamma}|\la\pi_\lambda(s)\pi_\lambda(s).v, \pi_\lambda(\gamma).w\ra|
=\sum_{\gamma\in \Gamma}|f_{\pi_\lambda(s).v}(\gamma^{-1}s )|\cr 
&\leq C_Q \|f_{\pi_\lambda(s).v}\|_1\leq C_Q C_{t\over 2} \|v\|.\cr}$$ 
This proves the proposition.\qed

Now we are going to show that the Poincar\'e series of the highest
weight vector 
$P(v_\lambda)$ is non-zero for almost all parameters $\lambda\in 2\Lambda$. 
For that we have to recall some facts concernig the generalized Harish-Chandra 
decomposition. Our source of reference is [Ne99a, Sect.\  XII.1]. 

We denote by $K$ the analytic subgroup of $G$ corresponding to $\k$. Note that 
$K$ is compact by our assumptions on $G$. Let $K_\C$ denote the universal 
complexification of $K$ and note that $K\subeq K_\C$. By $G_\C$ we denote the universal 
complexification of $G$. 

\par We define subalgebras of $\g_\C$ by 
$$\p^+\:=\oplus_{\alpha\in \Delta_n^+} \g_\C^\alpha\quad \hbox{and}\quad \p^-\:=\oplus_{\alpha\in
\Delta_n^-}\g_\C^\alpha$$   
and note that $\p^\pm$ are abelian since $W$ was assumed to be pointed. 
Further we have $\g_\C=\p^-\oplus\k_\C\oplus \p^+$. 

\par Let $P^\pm$ denote the analytic subgroups of $G_\C$ corresponding to 
$\p^\pm$. Assume for the moment 
that $G\subeq G_\C$. Then the multiplication mapping 
$$ P^- \times K_{\C} \times P^+ \to P^-K_\C P^+, \ \ (p_-, k, p_+)\mapsto p_-kp_+. $$
is biholomorphic and we have $G\subeq \oline S\into  P^-K_\C P^+$. 
\par For the general case we can use some standard 
covering theory to lift the results from above: We obtain a complex manifold 
$P^- K_\C P^+$ with a biholomorphic map $ P^- \times K_{\C} \times P^+ \to P^-K_\C P^+$
and the inclusion chain from before also lifts: $G\subeq \oline S\into  P^-K_\C P^+$. 
\par We denote by $\kappa\: \oline S\to K_\C$ the middle projection. 
Set $\t_0=\t\cap[\k,\k]$ and note that $\t=\t_0\oplus\z(\k)$. Accordingly we have 
$\t^*=\t_0^*\oplus \z(\k)^*$.

\Proposition 3.9. Let $\lambda_0\in i\t_0^*$ be dominant integral for 
$\Delta_k^+$. 

\item{(i)} There exists a lattice $\Gamma_{\z(\k)}\subeq i\z(\k)^*$ such that 
$$\lambda_0 +(\Gamma_{\z(\k)}\cap i C^\star)\bs\{0\}\subeq 2\Lambda.$$

\item{(ii)} Fix $\zeta\in HW(G,\Delta^+)\cap \big((\Gamma_{\z(\k)}
\cap i C^\star)\bs \{0\}\big)$ and set $\lambda_n\:=
\lambda_0 +n\zeta$. Then there exists an $N_0\in \N$ such that 
$\lambda_n\in  2\Lambda$ for all $n\geq N_0$, $n\in\N$. 
Further there exists an $n\in \N$, $n\geq N_0$ such that 
$$\la P(v_{\lambda_n}), v_{\lambda_n}\ra \neq 0, $$ 
where $v_{\lambda_n}$ is a normalized highest weight vector for $(\pi_{\lambda_n}, {\cal H}_{\lambda_n})$.

\Proof. (i) This follows from the structure of the set of unitary 
highest weights (cf.\ [Kr99, Sect.\ IV]) together with our explicit 
description of $\Lambda$.

\par\nin (ii) From the definition of $\Lambda$ it is clear that there 
exists an $N_0\in \N$ such that $\lambda_n\in 2\Lambda$ for all 
$n\geq N_0$, $n\in \N$. Define $\N_{\geq N_0}\:=\N\cap [N_0, \infty[$. 
We consider the function 
$$F\: \N_{\geq N_0}\to \C, \ \ n\mapsto \la P(v_{\lambda_n}), v_{\lambda_n}\ra=
\sum_{\gamma\in\Gamma}\la \pi_{\lambda_n}(\gamma).v_{\lambda_n}, v_{\lambda_n}\ra.$$ 
Note that the defining series converges absolutely by Proposition 3.8. 
\par If $(\pi_\lambda, {\cal H}_\lambda)$ is any unitary highest weight representation of $G$, 
then $F(\lambda)\:=\span\{\pi(K).v_\lambda\}$ is a finite dimensional 
irreducible representation of $K$. Hence the representation  $(\pi_\lambda, F(\lambda))$ of $K$ 
naturally extends to a holomorphic representation of $K_\C$. 
By [Kr98b, Prop.\ II.20] we now have 
$\langle \pi_{\lambda} (s). v_{\lambda}, 
v_{\lambda} \rangle =\langle \pi_{\lambda} (\kappa(s)).v_{\lambda}, 
v_{\lambda} \rangle$ for all $\lambda\in HW(G,\Delta^+)$ and $s\in \oline S$, and so  
$$(\forall n\in \N) \qquad F(n)=
\sum_{\gamma\in \Gamma} \la \pi_{\lambda_n}(\kappa(\gamma)).v_{\lambda_n}, v_{\lambda_n}\ra.
\leqno(3.3)$$ 

Let $\n_k := \bigoplus_{\alpha \in \Delta_k^+} \g_{\C}^{\alpha}$, 
$\overline{\n_k} := \bigoplus_{\alpha \in \Delta_k^-} \g_{\C}^{\alpha}$
be the complex conjugate and $N_k$, $\overline{N_k}$, the
respective analytic subgroups of $K_{\C}$. Then we have the
Bruhat decomposition
$$ K_{\C} = \bigcup_{w \in {\cal W}_{\k}} \overline{N_k}
\tilde{w} T_\C N_k\, ,$$ 
where $\tilde{w} \in \tilde K_{\C}$ represents $w$. 
Fix $n\in\N$ and consider the function 
$$f\: K_\C\to \C, \ \ k\mapsto \langle \pi_{\lambda_n} (k). v_{\lambda_n}, 
v_{\lambda_n} \rangle.$$
The fact that $\lambda_0$ is regular implies that $f$ vanishes 
on all cells $\overline{N_k}\tilde{w} T_\C N_k$ with $\tilde w\neq \1$. 
Write an element $k\in \overline{N_k} T_\C N_k$ as $k=\oline n (k) t(k) n(k)$
and note that $f(k)=f(t(k))=t(k)^{\lambda_n}$. Thus we have 
$$(\forall k\in K_\C)\quad f(k)=\cases{ 0 & if $k\not\in  \overline{N_k} T_\C N_k$,\cr
t(k)^{\lambda_n} & if $k\in  \overline{N_k} T_\C N_k$\cr}.\leqno(3.4)$$
\par The set $\{t(\kappa(\gamma))\in T_\C\:\gamma\in \Gamma\}$ defines 
a countable family $(t_j)_{j\in \N}\subeq T_\C$. Then (3.3) and (3.4) give that 
$$F\: \N_{\geq N_0}\to\C, \ \ n\mapsto \sum_{j=1}^\infty t_j^{n(\lambda_0 +\zeta)}.$$
Now it follows from Lemma 3.10 below that $F(n)\neq 0$ for at least one $n\in\N$
proving the proposition. \qed

\Lemma 3.10. Let $V$ be a finite dimensional real vector space, $V^*$ its dual and 
$V_\C=V+iV$ its complexification. Let $C\subeq V^*$ be a convex cone and 
$\Gamma\subeq V^*$ be a lattice. Suppose that $C\cap\Gamma\neq \eset$ and that 
there exists a sequence 
$(z_n)_{n\in \N}$ in $V_\C$ such that 
$$F\: C\cap \Gamma\to\C, \ \ \lambda\mapsto\sum_{n\in\N} e^{\lambda(z_n)}\, ,$$
is defined by absolutely  convergent series. Then $F\neq 0$. 

\Proof. The assertion of the lemma easily reduces to the one dimensional 
case, i.e., $V=\R$. By rescaling we may assume that 
$\Gamma=\Z$. Then $C=\R$ or $C=\pm]0,\infty[$ and taking a subcone of $C$ we,
up to sign,  may 
assume that  $C=]0,\infty[$. Then we have $\Z\cap ]0,\infty[=\N$ and  
$$F\: \N\to\C, \ \ m\mapsto \sum_{n=1}^\infty e^{mz_n}\, ,$$
with $\sum_{n=1}^\infty e^{m\Re(z_n)}<\infty$ for all $m\in \N$. 
Hence $\Re(z_n)>0$ for only finitely many $n\in \N$ and 
$$x_0\:=\sup_{n\in\N} \{ \Re z_n\: n\in\N\}<\infty.$$
Further $F$ extends to a bounded analytic function 
$$f\: ]1,\infty[\to\C, \ \ \lambda\mapsto \sum_{n=1}^\infty e^{\lambda z_n}.$$
To obtain a contradiction, assume that $F=0$. Then $f=0$ by [Kr01, App.\ A]. 

\par W.l.o.g.\ we may assume that $\Re z_n=x_0$ exactly for 
$1\leq n\leq N$ for some $N\in\N$ and that 
$z_1=z_j$ exactly for $1\leq j\leq k$, $k\leq N$. 
Then dominated convergence gives that 
$$\lim_{\lambda\to\infty} e^{-\lambda z_1}\Big(\sum_{j=N+1}^\infty e^{\lambda z_j}\Big)=0.$$ 
Choose $\lambda_0>0$ such that 
$|e^{-\lambda z_1}\Big(\sum_{j=N+1}^\infty e^{\lambda z_j}\Big)|<{1\over 2}$ for all $\lambda>\lambda_0$. 
Then $f=0$ implies that 
$$(\forall \lambda>\lambda_0)\quad |k+\sum_{j=k+1}^N e^{i\lambda \Im (z_j-z_1)}|<{1\over 2}.$$
Since the one parameter subsemigroup  
$$\{ (e^{i\lambda\Im (z_{k+1}-z_1)}, \ldots, e^{i\lambda\Im (z_N-z_1)})\: \lambda>\lambda_0\}$$
in $(\SS^1)^{N-k}$ contains elements arbitary close to the identity  $\1\in (\SS^1)^{N-k}$,  
we arrive at a contradiction, proving the lemma.\qed

We now summarize the main results of this section.

\Theorem 3.11. Let $S=\Gamma_G(W)$ be a complex Ol'shanski\u\i{} semigroup with 
$G$ a (CA)-Lie group and $Z(G)$ compact. Let $U$ be any open 
cone in $iC^\star$. Then for every discrete 
subgroup $\Gamma<G$ there exist a $\Gamma$-spherical contractive
unitary  highest weight representation 
$(\pi_\lambda, {\cal H}_\lambda)$ of $G$ with $\lambda\in U$ and compact kernel.

\Proof. Note that non-zero $\Gamma$-fixed hyperfunction vectors for a
unitary highest weight representation $(\pi_\lambda, {\cal H}_\lambda)$ are 
cyclic, since ${\cal H}_\lambda^{T,\omega}$ is dense in ${\cal H}_\lambda^\omega$ 
(cf.\ [HiKr99, Lemma 6.2.1]). Therefore it follows from Proposition 3.8 and Proposition 3.9 (ii) that there 
exists a $\Gamma$-spherical unitary highest weight representation of $G$. Moreover 
the representations belonging to a parameter $\lambda\in 2\Lambda$ are 
contractive since $\Lambda\subeq iC^\star$ (cf.\ [Kr99, Lemma IV.13]).  
The remaining statement of the theorem now follows from Proposition 3.9(ii)
together with [Ne99a, Lemma X.4.7, Th.\ X.4.8].\qed

\sectionheadline {4. Point separation and vansihing at infinity}

In this section we construct rich classes of holomorphic functions on 
$\Gamma\bs S$. For $\Gamma<G$ cocompact we construct 
a class of point separating holomorphic functions on 
$\Gamma\bs S$ which vanish at infinity. For arbitray $\Gamma$
we will show that $\Hol(\Gamma\bs S)$ separates the points.

\subheadline{Compact $G$-orbits in the hyperfunctions}

\Lemma 4.1. Let $s_0=\Exp(X_0)$ with $X_0\in iC$. Then the 
following assertions hold:

\item{(i)} $Ss_0\subeq G\Exp\big(\conv({\cal W}_\k.X_0) +iC\big)G$.
\item{(ii)} There exists an $X\in iC$ such that 
$$Ss_0\subeq G\Exp(X+iC)G.$$

\Proof. (i) For $Z\in iC$ set $F_{Z,X_0}\:=\conv({\cal W}_\k.(Z+X_0))+ iC $. Then 
[Ne97, Prop.\ V.1] shows that 
$$(\forall Z\in i C)\quad \Exp(Z)G\Exp(X_0)\subeq G\Exp(F_{Z,X_0})G.$$
Now we have 
$${\cal W}_\k.(Z+X_0)\subeq {\cal W}_\k.Z+{\cal W}_\k.X_0\subeq {\cal W}_\k.X_0+ iC
\subeq \conv({\cal W}_\k.X_0)+iC$$
for all $Z\in iC$. Hence  (i) follows from $S=G\Exp(iC) G$.

\par\nin (ii) Let $\leq_C$ denote the conal order on the convex cone $iC$, i.e., 
for $X,Y\in iC$ one has $X\leq_C Y$ if and only if $X\in Y+ iC$. 
Since 
$\conv({\cal W}_\k.X_0)$ is compact, we therefore find an element $Y\in iC$ such that 
$Y\leq_C \conv({\cal W}_\k.X_0)$, i.e.,  $\conv({\cal W}_\k.X_0)\subeq Y+iC$. In view of 
(i), this proves (ii).\qed

\Lemma 4.2. Let $(\pi,{\cal H})$ be a holomorphic representation of $S$. 
Suppose that $\Gamma\bs G$ is compact and let $v^\Gamma\in ({\cal H}^{-\omega})^\Gamma$. 
Then the following assertions hold: 

\item{(i)} The orbit $\pi^{-\omega}(G).v^\Gamma$ is compact in ${\cal H}^{-\omega}$.
\item{(ii)} For all $X_0\in iC$ the set
$\pi^{-\omega}\big(\Exp(X_0)\big)\pi^{-\omega}(G).v^\Gamma$ is compact in ${\cal H}$.

\Proof. (i) Recall from Definition 3.1 that the representation
$(\pi^{-\omega}, {\cal H}^{-\omega})$ of $G$ is continuous, i.e., the map 
$$G\times {\cal H}^{-\omega}\to {\cal H}^{-\omega},\ \ (g,v)\mapsto 
\pi^{-\omega}(g).v \, ,$$
is continuous. Thus $\pi^{-\omega}(G).v^\Gamma$ is compact by the cocompactness of $\Gamma$ in 
$G$.

\par\nin (ii) Recall the definition of the topology on ${\cal H}^{\omega}$
(cf.\ Remark 3.2). Being the strong dual of the analytic vectors 
${\cal H}^\omega$, the topology on ${\cal H}^{-\omega}$ is the finest 
locally convex topology for which all maps 
$${\cal H}^{-\omega}\to {\cal H}, \ \ v\mapsto \pi^{-\omega}(\Exp(tX_0)).v,
\quad t>0,$$
become continuous. Hence for all compact subsets 
$K\subeq {\cal H}^{-\omega}$, the set $\pi^{-\omega}(\Exp(X_0)).K$ is 
compact in ${\cal H}$. In view of (i), this proves (ii).\qed

\subheadline{Realizing $\Gamma$-spherical representations in holomorphic functions}

Throughout this subsection $(\pi,{\cal H})$ denotes a holomorphic representation of 
a complex Ol'shanski\u\i{} semigroup. We assume that $(\pi,{\cal H})$ is 
$\Gamma$-spherical for an arbitrary discrete subgroup $\Gamma<G$. 

Fix now a cyclic element 
$v^\Gamma\in ({\cal H}^{-\omega})^\Gamma$. 
By [KN\'O97, App.] the mapping 
$$r\: {\cal H}\to \Hol(\Gamma\bs S), \ \
 v\mapsto f_v;\  f_v(\Gamma s)\:=\la \pi(s).v,v^\Gamma\ra\, , \leqno(4.1)$$
is injective, continuous and $G$-equivariant, i.e., the 
map $r$ gives us a realization of ${\cal H}$ as a $G$-invariant 
Hilbert space of holomorphic functions on $\Gamma\bs S$ (cf.\ Definition 2.5). In the 
sequel we assume that ${\cal H}\subeq \Hol(\Gamma\bs S)$. 
Recall from Section 2  that there is a reproducing kernel $K$ for ${\cal H}$.

\Lemma 4.3. Let $(\pi, {\cal H})$ be a $\Gamma$-spherical holomorphic
representation of $S$ realized in the holomorphic functions on $\Gamma\bs S$ 
via 
the realization {\rm (4.1)}. If $K$ is the reproducing kernel of ${\cal H}$, then  
$K_{\Gamma s}=\pi(s^*).v^\Gamma$, $s\in S$, and 
$$(\forall s,t\in S)\qquad  K(\Gamma s,\Gamma t)=\la \pi(t^*).v^\Gamma, \pi(s^*).v^\Gamma\ra.$$ 

\Proof. This is analogous to the proof of [HiKr01, Th.\ 4.1.1]. \qed

\Proposition 4.4. Assume that $\Gamma\bs G$ is compact and let $(\pi, {\cal H})$ be 
a $\Gamma$-spherical holomorphic contraction representation of $S$. 
Then the following assertions hold: 

\item{(i)} Fix $X_0\in iC$ and set $s_0\:=\Exp(X_0)$. Then, the reproducing kernel $K$
of ${\cal H}$ is bounded when restricted to $\Gamma\bs Ss_0\times \Gamma\bs Ss_0$. 
In particular, there exists a constant $C=C(s_0)$ such that 
$$(\forall s\in Ss_0)\qquad  \| K_{\Gamma s}\|\leq C.$$

\item{(ii)} All analytic vectors ${\cal H}^\omega$ are bounded holomorphic 
functions on $\Gamma\bs S$.

\Proof. (i) By Lemma 4.3 we have $K(\Gamma s,\Gamma t)=\la \pi(t^*).v^\Gamma, \pi(s^*).v^\Gamma\ra.$ 
Since $|K(z,w)|^2\leq K(z,z) K(w,w)$ for all $z,w\in \Gamma\bs S$ (cf.\ [Ne99a, Ch.\ I]), we have 
$$\sup_{z,w\in \Gamma\bs S} |K(z,w)|\leq \sup_{z\in \Gamma\bs S} K(z,z).$$
In view of  Lemma 4.1 we have 
$Ss_0\subeq G\Exp(X)\Exp(iC)G$ for some $X\in iC$. Therefore we get for all 
$s\in S$  that 
$$\|K_{\Gamma ss_0}\|^2=K(\Gamma s s_0, \Gamma s s_0)=\|\pi(ss_0).v^\Gamma\|^2\leq 
\sup_{g\in G}\|\pi(\Exp(X))\pi(g).v^\Gamma\|^2 ,$$  
since $(\pi, {\cal H})$ is contractive and $\Exp(X)$ commutes with $\Exp(iC)$. 
But now the assertion follows from Lemma 4.2(ii). 

\par\nin (ii) Let $f\in {\cal H}^\omega$. Then $f\in \pi(s_0).{\cal H}$ 
for some 
$s_0=\Exp(X_0)$, $X_0\in iC$. Thus we can write $f=\pi(s_0).g$ with $g\in {\cal H}$. 
But then (i) implies that 
$$\sup_{s\in S} |f(\Gamma s)|=\sup_{s\in Ss_0} |g(s)|=
\sup_{s\in Ss_0}|\la g, K_{\Gamma s}\ra|\leq C(s_0)\|g\|.\qeddis

\subheadline{Vanishing at infinity}

In this subsection we show that the analytic vectors of certain
invariant Hilbert spaces of holomorphic functions on $\Gamma\bs S$
are vanishing at infinity. 
 
\Lemma 4.5. Let $s_0=\Exp(X_0)\in S$, $X_0\in iC$, and 
$(s_n)_{n\in \N}$ be a sequence in $Ss_0$. 
\item{(i)} There exists an element $X\in iC$ such that 
$$s_n=g_nh_n^{-1}\Exp(X+X_n)h_n$$
with $g_n, h_n\in G$, $X_n\in iC$.
\item{(ii)} Assume now that $\Gamma< G$ is cocompact and  that 
$(\Gamma s_n)_{n\in\N}$ leaves every compact subset in 
$\Gamma\bs \oline S$. Then the sequence $h_n^{-1}\Exp(X+X_n)h_n$
leaves every compact subset of $\Exp(\oline W)$. In particular 
we have $X_n\to\infty $ in $i\oline C$ or $h_n\to\infty$ in G.

\Proof. (i) This follows from Lemma 4.1. 
\par\nin (ii) Note that there is a polar decomposition of $\Gamma\bs \oline S$
$$\Gamma\bs G \times \oline W\to \Gamma\bs \oline S$$ 
which is a homeomorphic mapping. As $\Gamma\bs G$ is compact, (ii) follows now from (i). \qed

\Theorem 4.6. {\rm (Vanishing at Infinity)} Let $\Gamma< G$ be cocompact 
and $(\pi_\lambda, {\cal H}_\lambda)$ be a $\Gamma$-spherical unitary 
highest weight representation of $G$ with $\lambda\in i\Int C^\star$. 
Assume further that $G$ is a (CA)-Lie group and that $Z(G)$ is compact. 
Let $v\in {\cal H}^\omega$ and 
$f_v(\Gamma s)=\la \pi(s).v,v^\Gamma\ra$ the corresponding holomorphic function on 
$\Gamma\bs S$. Then $f_v$ extends to a continuous function on $\Gamma\bs\oline S$, also 
denoted by $f_v$ and the 
following assertions hold:

\item{(i)} The function $f_v$ is bounded. 
\item{(ii)} We have $\lim_{s\to\infty\atop s\in \Gamma\bs \oline S} f_v(s)=0$. 

\Proof. Fix $X_0\in iC$ and set $s_t\:=\Exp(tX_0)$ for $t>0$. Then ${\cal H}_\lambda^\omega =\cup_{t>0} 
\pi_\lambda(s_t).{\cal H}_\lambda$ shows that there is a $t_0>0$ such that $v=\pi_\lambda(s_{t_0}).w$ for some $w\in {\cal H}_\lambda$. 
Set $s_0\:=s_{t_0}$. Then $\oline Ss_0\subeq S$ shows that there is a continuous extension of $f_v$ to 
$\Gamma\bs \oline S$. 
\par\nin (i) Proposition 4.4(ii).
\par\nin (ii) Let $(s_n)_{n\in \N}$ be a sequence in $S$ such that 
$\Gamma s_n\to \infty$ in $\Gamma\bs \oline S$. 
We claim that $\Gamma s_ns_0 \to \infty$ in $\Gamma\bs\oline S$. 
For that it is sufficient to show that the mapping 
$$\rho_\Gamma(s_0)\: \Gamma\bs \oline S\to \Gamma\bs \oline, 
\ \ \Gamma s\mapsto \Gamma s_0$$
is proper. Recall from [HiNe93, Th.\ 3.20] that 
$\oline S$ acts via 
right multiplication on $\oline S$ by proper maps. 
In particular, $\rho(s_0)\: \oline S\to \oline S, \ s\mapsto ss_0$
is proper. Let now $Q_\Gamma\subeq \Gamma\bs \oline S$
be a compact subset. Then we find a compact subset $Q\subeq \oline S$
such that $Q_\Gamma=\Gamma\bs \Gamma Q$. Since $\rho(s_0)$ is proper, 
the set $\rho(s_0)^{-1}(Q)$ is compact. Hence 
$\rho_\Gamma(s_0)^{-1}(Q_\Gamma)=\Gamma\bs \Gamma\rho(s_0)^{-1}(Q)$
is compact, establishing our claim. 

\par By Lemma 4.5 we find an element $X\in C$ such that 
$$s_n s_0 =g_nh_n^{-1}\Exp(X+X_n)h_n$$
for elements $g_n, h_n\in G$, $X_n\in iC$. 
Then Lemma 4.3 and Proposition 4.4 show that  
$$\eqalign{f_v(\Gamma s_n)&=f_w(\Gamma s_n s_0)=\la w, K_{\Gamma s_n s_0}\ra =
\la w, \pi_\lambda(s_0s_n^*).v^\Gamma\ra\cr 
& =\la w, \pi_\lambda(h_n^{-1}\Exp(X+X_n))\pi_\lambda^{-\omega}(h_ng_n^{-1}).v^\Gamma\ra\cr
&=\la \pi_\lambda(h_n).w, \pi_\lambda(\Exp({1\over 2}X+X_n))\pi_\lambda^{-\omega}(\Exp({1\over 2} X)
h_ng_n^{-1}).v^\Gamma\ra.\cr}$$
Since $\Gamma\bs G$ is compact, Lemma 4.2(ii) implies that $Q\:=\pi_\lambda^{-\omega}
(\Exp({1\over 2} X)G).v^\Gamma$
is compact in ${\cal H}_\lambda$. 
Thus we get 
$$|f_v (\Gamma s_n)|\leq \sup_{u\in Q} |\la \pi_\lambda(h_n).w, \pi_\lambda(\Exp({1\over 2}X+X_n)).u\ra|.
\leqno(4.2)$$
As $\Gamma s_n$ tends to infinity in $\Gamma\bs \oline S$ we have $h_n\to\infty $ in $G$ or 
$X_n\to \infty$ in $\oline C$ (cf.\ Lemma 4.5(ii)). 

\par\nin Case 1: $X_n\to\infty$. Then ${1\over 2} X+ X_n\to\infty $ in $i\oline C$ and we have  
$$\|\pi_\lambda(\Exp({1\over 2}X +X_n))\|=\sup_{w\in {\cal W}_\k} e^{w.\lambda({1\over 2}X +X_n)}$$
(cf.\ [Kr99, Lemma 4.13]). Thus $\lambda \in i\Int C^\star$ and 
the fact that $C$ 
is ${\cal W}_\k$-invariant (and hence $i\Int C^\star$)
shows that $\|\pi_\lambda(\Exp({1\over 2}X +X_n))\|\to 0$. 
It follows from (4.2) that 
$$\lim_{n\to \infty}|f_v(\Gamma s_n)|\leq \|w\|\cdot\sup_{u\in Q}\|u\|\cdot\lim_{n\to\infty} 
\|\pi_\lambda(\Exp({1\over 2}X +X_n))\| =0.$$
\par\nin Case 2: $h_n\to\infty$. 
We may assume that $X_n$ stays in a compact subset of $i\oline C$ because of Case 1. 
Then $\pi_\lambda(\Exp({1\over 2}X+X_n)).Q$ stays in a compact subset $Q'\supeq Q$ of ${\cal H}_\lambda$ 
and (4.2) implies that 
$$|f_v(\Gamma s_n)|\leq \sup_{u\in Q'} |\la \pi_\lambda(h_n).w, u\ra|.$$
Note that $(\pi_\lambda, {\cal H}_\lambda)$ has compact 
kernel, since $G$ is a (CA)Lie group, $Z(G)$ is compact and
$\lambda\in i\Int C^\star$ (cf.\ the proof of [Kr98a, Prop.\ V.7]). 
Hence the facts that  $(\pi_\lambda, {\cal H}_\lambda)$ is irreducible 
and $G$ is a (CA)-Lie 
group imply that Mayer's generalization of the 
Howe-Moore Theorem (cf.\ [Ma97, Prop.\ 3.4], 
[HoMo79, Th.\ 5.1]) applies and so $\lim_{n\to\infty}|f(\Gamma s_n)|=0$. 
\qed

\subheadline{Separating holomorphic functions}

Now we are ready to prove the main result in this section. 
Recall from [Ko98, Ch.\ 3, \S2] the definition of 
the Caratheodory semimetric and note that the Caratheodory semimetric 
is a metric if and only if  the bounded 
holomorphic functions separate the points.

\Theorem 4.7. {\rm (Separation of Points)} 
Let $S$ be a complex Ol'shan\-ski\u\i{} semigroup associated to 
a (CA)-Lie group $G$. Let $\Gamma<G$ be an arbitrary discrete subgroup of 
$G$. Then the following assertions hold: 

\item{(i)} For every $z_1, z_2\in \Gamma\bs S$, $z_1\neq z_2$, there exists an $f\in {\cal B}^1(\Gamma\bs S)$
such that $0\neq f(z_1)\neq f(z_2)$. In particular $\Hol(\Gamma\bs S)$ separates the points of $\Gamma\bs S$. 
\item{(ii)} If $\Gamma<G$ is cocompact, then the bounded holomorphic functions on $\Gamma\bs S$ separate 
the points. In particular, the Caratheodory semimetric on $\Gamma\bs S$ is 
a metric.

\Proof. (i) The proof uses a trick which was used by Siegel in order to show that 
the Siegel modular forms separate the points on $G/K$ (cf.\ [Fr83]). 
\par   Let $C_0(\Gamma\bs \oline S)$ denote the algebra of continuous 
function on $\Gamma\bs \oline S$ which vanish at infinity. Then it follows from 
[Kr98a, Prop.\ II.4(i)(b)] and a slight adapation of [Kr98a, Lemma V.13] that 
${\cal B}^1(S)\cap C_0(\oline S)$ separates the points of $S$. 

\par Let now $\Gamma s_1, \Gamma s_2\in \Gamma\bs S$ be two different points. Let 
$f\in {\cal B}^1(S)\cap C_0(\oline S)$ be such that 

$$m\:=\sup_{\gamma\in \Gamma} \{ |f(\gamma s_1)|, |f(\gamma s_2)|\}>0.$$
By Proposition 3.4 we know that for every $f\in {\cal B}^1(S)$ and $s\in S$ one has 

$$|f^\Gamma(s)|\leq \sum_{\gamma\in \Gamma} |f(\gamma s)|<\infty\leqno(4.3)$$
and $f^\Gamma\in {\cal B}^1(\Gamma\bs S)$. 

Hence the supremum $m$ is actually attained and without loss of generality we may assume that 
$m=1=f(s_1)$. 
Let now 
$$\Gamma_0\:=\{ \gamma\in \Gamma\bs \{\1\} \: |f(\gamma s_1)|=1 \quad\hbox{or}\quad   |f(\gamma s_2)|=1\} $$
and note that $\Gamma_0$ is a finite set by (4.3). 
Since ${\cal B}^1(S)\cap C_0(\oline S)$ separates the points of $S$ we find a 
$g\in {\cal B}^1(S)\cap C_0(\oline S)$ such that 
$$g(s_1)=1 \quad \hbox {and} \quad g(\gamma s_1)=g(\gamma s_2)=0\quad  (\forall\gamma\in \Gamma_0).$$ 
For every $n\in \N$ we now consider the function 
$$F_n(s)=g(s) f^n(s)$$
and note that $F_n\in {\cal B}^1(S)$. Averaging over $\Gamma$ we obtain that 

$$(F_n)^\Gamma(\Gamma s) =\sum_{\gamma\in \Gamma} g(\gamma s) f^n(\gamma s)$$
for all $s\in S$. Note that $(F_n)^\Gamma\in {\cal B}^1(\Gamma\bs S)$. Now we have 

$$|(F_n)^\Gamma|(\Gamma s_2)\leq \sum_{\gamma\in \Gamma\bs \Gamma_0} 
|g(\gamma s_2)||f^n(\gamma s_2)|$$
and we see that $\lim_{n\to\infty}  (F_n)^\Gamma(\Gamma s_2)=0$ by 
(4.3) and $|f(\gamma s_2)|<1$ for all $\gamma\in \Gamma\bs \Gamma_0$. 
On the other hand we have

$$(F_n)^\Gamma(\Gamma s_1)= g(s_1) f(s_1) +\sum_{\gamma\in \Gamma\bs \Gamma_0} 
g(\gamma s_2)f^n(\gamma s_2)= 1 +\sum_{\gamma\in \Gamma\bs \Gamma_0} 
g(\gamma s_2)f^n(\gamma s_2)$$
and by the same reason as before we have $\lim_{n\to\infty}  (F_n)^\Gamma(\Gamma s_1)=1$. 
Hence there exists an $n\in \N$ such that $(F_n)^\Gamma(\Gamma s_1)\neq (F_n)^\Gamma(\Gamma s_2)$
proving (i) 
\par\nin (ii) In view of the proof of (i), this follows from Proposition 4.4. \qed  

\Remark 4.8. By Neeb's version of the Gelfand-Raikov theorem 
(cf.\ [Ne99, Th.\ XI.5.1]) 
for complex Ol'shanski\u\i{} semigroups, 
the bounded holomorphic functions on $S$ 
separate the points. It follows that $S$ is hyperbolic complex. 
Hence $\Gamma\bs S$ is hyperbolic complex for arbitrary discrete subgroups 
$\Gamma<G$ by [Sh92, \S 20, Th.\ 1]. To compare this
with the results above  notice that Theorem 4.7(ii) is a stronger 
statement than that $\Gamma\bs S$ is hyperbolic complex. \qed

\Example 4.9. (a) Let $G$ be a linear Hermitian Lie group. 
Then by a theorem of 
Borel (cf.\ [Bo63, Th.\ C], [Ra72, Th.\ 14.1]) 
there exist cocompact lattices $\Gamma<G$. Further the Lie algebra 
of $\g$ admits pointed open convex $\Ad(G)$-invariant elliptic cones $\eset\neq W\subeq \g$. 
Set $S=\Gamma_G(W)$. Then by Theorem 4.8 the bounded holomorphic functions 
on $\Gamma\bs S$ separate points. 
 
\par\nin (b) Let $\h_1=\R X\oplus \R Y\oplus \R Z$ be the 3-dimensional 
Heisenberg algebra  with relations $[X,Y]=Z$ and all other brackets 
vanishing. Then $\a=\R T$ acts on $\h$ by derivations 
via $[T,X]=Y$, $[T,Y]=-X$ and $[T,Z]=0$. 
Then $\g=\h\rtimes\a$ is a four dimensional solvable Lie algebra, the so-called 
{\it oscillator algebra}. It is the basic example of a solvable Lie 
algebra admitting  pointed open invariant convex cones $\eset\neq W\subeq\g$. 
Note that $\t=\R Z\oplus\a$ is a compactly embedded Cartan subalgebra of $\g$. 
\par Let $H_1\cong\R^3$ denote the simply connected Heisenberg group corresponding to $\h_1$. 
We have a connected Lie group $G=H_1\rtimes \SS^1$ with Lie algebra $\g$. Then 
$\Gamma=\Z^2\times {1\over 2}\Z\subeq H_1$ is a cocompact lattice in $G$ and hence the 
bounded holomorphic functions 
on $\Gamma\bs S$ separate points by Theorem 4.8.

\par\nin (c) Now we give an example of a Lie group $G$ of mixed type, i.e., $G$ is neither 
reductive nor solvable. We define a quadratic form 
$$F\:\R^3\to\R, \ \ (x,y,z)\mapsto x^2+y^2-3z^2$$
and let $L\:=\SO_0(F,\R)$ denote the connected component of the special invariance 
group of $F$. Then $F$ does not represent zero in a non-trivial 
rational way so that $\Gamma_L\:=\SO_0(F,\Z)$ forms a cocompact lattice in 
$L$ (cf.\ [Bo69, 8.6, Ex.\ 1]). 
\par Let $\la\cdot, \cdot\ra_{2,1}$ denote the non degenerate symmetric 
bilinear form on $\R^3$ induced by $F$. Then the prescription 
$$\Omega_F\: (\R^3\times\R^3)\times (\R^3\times\R^3)\to\R, \ \ \big((u,v), (u',v')\big)\mapsto 
\la u,v'\ra_{2,1}- \la u',v\ra_{2,1}$$
defines an $L$-invariant symplectic form on $\R^3\times\R^3$ (we let $L$ act diagonally). 
Associated to $\Omega_F$ we build the nilpotent Lie algebra 
$\uu\:=(\R^3\times \R^3)\times \R$ with bracket 
$$[((u,v), z), ((u',v'), z')]=\big(0, \Omega_F((u,v), (u',v'))\big).$$ 
Note that $\uu$ is isomorphic to the 7-dimensional Heisenberg algebra $\h_3$. Let 
$U$ denote a simply connected Lie group with Lie algebra $\uu$. Then 
$U\cong (\R^3\times\R^3)\times\R$ and $\Gamma_U\:=(\Z^3\times\Z^3)\times {1\over 2}\Z$ is a
cocompact lattice in $U$ which is stable under $\Gamma_L$. Thus 
$\Gamma\:=\Gamma_U\rtimes \Gamma_L$ is a cocompact lattice in the 
semidirect product $G\:=U\rtimes L$. 
\par Note that there is an embedding of $\g$ into the Jacobi algebra 
$\h_3\rtimes \sp(3,\R)$ sending $\uu$ isomorphically onto $\h_3$ and $\l\:=\L(L)$ into $\sp(3,\R)$. 
In view of [Ne99a, Ch.\ VII], this 
embedding guarantees us the existence of non-trivial pointed $\Ad(G)$-invariant open convex cones
$W\subeq \g$. Now the bounded holomorphic functions 
on $\Gamma\bs \Gamma_G(W)$ separate points by Theorem 4.8. 
\qed

\sectionheadline{5. The Stein property of $\Gamma\bs S$ for  $\Gamma\bs G=\Sl(2,\Z)\bs \Sl(2,\R)$}

We let the group $G=\Sl(2,\R)$ act on the upper halfplane $\Pi^+\:=\{ z\in \C\: \Im z>0\}$
by means of M\"obius transformations: 

$$g.z={az+b\over cz+d}\qquad g=\pmatrix {a& b\cr c& d\cr}\in G, \ z\in \Pi^+.$$
Note that $G_\C =\Sl(2,\C)$. 
The compression semigroup of $\Pi^+$ 
$$\oline S\:=\{g\in G_\C\: g.\Pi^+\subeq \Pi^+\}$$
defines a closed complex Ol'shanski\u\i{} semigroup with interior 
$$S\:=\{g\in G_\C\: g.\oline {\Pi^+}\subeq \Pi^+\}.$$
Let $U=\pmatrix{0 & 1\cr -1&0\cr}$ be the infinitesimal generator of the maximal compact 
subgroup $K=\SO(2,\R)$ of $G$. Set $\R^-=]-\infty, 0[$. Then it is 
easily shown that 

$$S=G\exp(i\R^- U)G=G\exp(iW)$$
with $W=\Ad(G).\R^-U$ the lower light cone in $\g=\sL(2,\R)$. Note that $W$ and $-W$ 
are the only non-trivial $\Ad(G)$-invariant pointed convex cones in $\g$. 
\par Let $\Gamma=\Sl(2,\Z)$ be the modular group.
For $z\in \Pi^+$ set $q\:=e^{2\pi i z}$. Then the discriminant
$$\Delta\: \Pi^+\to\C, \ \ z\mapsto q\prod_{n=1}^\infty (1-q^n)^{24}$$
defines a modular form of weight $12$, i.e., $\Delta$ is holomorphic, vanishes at infinity 
and satisfies 

$$(c_\gamma z+d_\gamma)^{-12} \Delta(\gamma.z)=\Delta(z) \qquad z\in \Pi^+, \ \gamma=
\pmatrix{a_\gamma & b_\gamma\cr c_\gamma& d_\gamma\cr}.$$
For $g=\pmatrix{a& b\cr c& d\cr}\in G$ and $z\in \Pi^+$ we define 
$$\mu(g,z)\:=cz+d$$
and note that $\mu$ is a cocycle, i.e., $\mu(g_1g_2, z)=\mu(g_1, g_2.z)\mu(g_2,z)$
for all $g_1, g_2, \in G$, $z\in \Pi^+$. 
Then the description 

$$f_\Delta\: G\to \C, \ \ g\mapsto \mu(g,i)^{12}\Delta(g.i)$$
defines a left $\Gamma$-invariant analytic function on $G$ (cf.\ [Bo66]). 
Since $\oline S$ compresses $\Pi^+$, we see that $\mu$ extends to a continuous 
function $\mu(s,z)$, $s\in \oline S$, $z\in \Pi^+$, with no zeros. 
By the same reason the prescription $s\mapsto \Delta(s.i)$ defines a 
continuous function on $\oline S$. Hence we can analytically continue
$f_\Delta$ to a continuous function 
$$F_\Delta\: \oline S\to \C, \ \ s\mapsto \mu(s,i)^{12} \Delta(s.i)$$
which is holomorphic when restricted to $S$. Note that $F_\Delta$ is left 
$\Gamma$-invariant and hence factors to a function on $\Gamma\bs \oline S$, 
also denoted by $F_\Delta$. 
\par It is a consequence of [KO01, Th.\ 3.2] that $F_\Delta$ vanishes 
at infinity on $\Gamma\bs S$, i.e., 
$$\lim_{\Gamma s\to\infty\atop \Gamma s\in \Gamma\bs \oline S}  F_\Delta(\Gamma s)=0.\leqno(5.1)$$
Let us briefly scetch the proof of this fact. The discriminant 
$\Delta$ is bounded on a fundamental domain for $\Gamma$ in $\Pi^+$. Hence it
is sufficient to show that $\lim_{s\to \infty \atop s\in \oline S} \mu(s,i)=0$. 
For $G=\Sl(2,\R)$ this can be seen by direct calculation; more generally it 
follows from the fact that the middle projection $\kappa\: P^-K_\C P^+\to K_\C$ 
restricted to $\oline S$ is a proper map and has image $\kappa(\oline S)=K\exp(i\R^-U)$
(cf.\ [KO01, Prop.\ 1.2, Cor.\ 2.4]).   
We note that $F_\Delta$ has no zeroes on $\Gamma\bs \oline S$. Hence ${1\over |F_\Delta|}$
defines a continuous plurisubharmonic function on $\Gamma\bs \oline S$ and 
(5.1) implies that 
$$\lim_{\Gamma s\to\infty(in\  \Gamma\bs \oline S) \atop \Gamma s\in \Gamma\bs  S}  {1\over |F_\Delta|}
(\Gamma s)=\infty.\leqno(5.2)$$

\Theorem 5.1. Let $S$ be a complex Ol'shansk\u\i{} semigroup associated to $G=\Sl(2,\R)$ and let 
$\Gamma=\Sl(2,\Z)$. Then $\Gamma\bs S$ is Stein. 

\Proof. In view of Grauert's solution of the Levi problem (cf.\ Theorem 2.2), it is sufficient
to prove the existence of a strictly plurisubharmonic exhaustion function of $\Gamma\bs S$. 
\par Let $\psi$ be the non-negative biinvariant plurisubharmonic function on $S$ from 
Theorem 2.4. Then $\psi$ factorizes to a function on $\Gamma\bs S$ which we also 
denote by $\psi$. Then 
$$\phi\:=\psi+ {1\over |F_\Delta|}$$
defines a strictly plurisubharmonic function on $\Gamma\bs S$. 
To conclude the proof of the theorem it suffices to show that $\phi$ is  proper. 
Let $(z_n)_{n\in\N}$ be a sequence 
in $\Gamma\bs S$ leaving every compact subset of $\Gamma\bs S$. Then we either 
have $z_n\to z\in \big( \Gamma\bs \oline S \big) - (\Gamma\bs S)$ or 
$z_n\to \infty$ in $\Gamma\bs \oline S$. In the first case we have 
$\psi(z_n)\to\infty$ by Theorem 2.4 while in the latter case 
$\lim_{n\to\infty} {1\over |F_\Delta(z_n)|}=\infty$ by (5.2). 
This concludes the proof of the theorem. \qed

\Remark 5.2. (a) Note that by Theorem 2.3 now also
for all subgroups $\Gamma_0< \Sl(2,\Z)$ the manifold $\Gamma_0\bs S$ is Stein.
\par\nin (b) Let $G$ be a hermitian Lie group and $\Gamma<G$ a lattice. 
We claim that there exists no 
holomorphic cusp forms on $G/K$ without zeros if $\Gamma<G$ is cocompact or if the real rank ${\rm rank}(G)$ of $G$
is strictly larger than one. Since if $f$ were a holomorphic cusp form 
without zeros, then ${1\over f}$ would again be an automorphic form. 
If $\Gamma<G$ is cocompact, this 
is clear, and if $\Gamma<G$ is not cocompact and ${\rm rank}(G)>1$, then this is implied by the 
Koecher-principle (cf.\ [Fr83, Ch.\ I, Hilfssatz 3.5]). 
Now ${1\over f}$ would be of ``negative weight'' which is impossible since there are no 
automorphic forms of negative weight (if $\Gamma<G$ is cocompact, this is a consequence 
of the Plancherel Theorem for $L^2(\Gamma\bs G)$ and in the non-cocompact case this can be proved as in 
[Fr83, Ch.\ I, Satz 3.13]).
\par Therefore our construction for $\Sl(2,\R)$ cannot be generalized to arbitrary hermitian groups nor 
to the case of cocompact subgroups. 
\qed

\Conjecture 5.3.  We conjecture that $\Gamma\bs S$ is Stein for all complex 
Ol'shanski\u\i{} semigroups and all discrete subgroups $\Gamma<G$. 
In view of the proof of Theorem 5.1, it would be sufficient to find a holomorphic function
on $\Gamma\bs S$ with no zeros and which vanishes at infinity on $\Gamma\bs \oline S$. 
If $\Gamma<G$ is cocompact, then we have seen in Theorem 4.6 that all matrix coefficients 
$f_v(\Gamma s)=\la \pi(s).v, v^\Gamma\ra$ for $v\in {\cal H}^\omega$ 
vanish at infinity. The difficulty is to show that there exists a certain $v\in {\cal H}^\omega$
so that $f_v$ has no zeros. As we pointed out in Remark 5.2, this can never happen 
if $v$ is a highest weight vector. But it is likely to be true that $f_v$ has no zeros 
for all $v\in \pi_\lambda(S).v^\Gamma$, as the reproducing kernel 

$$K(\Gamma s, \Gamma t)=\la \pi(t^*).v^\Gamma,\pi(s^*).v^\Gamma\ra$$    
of the realiztion of $(\pi, {\cal H})$ in $\Hol(\Gamma\bs S)$ should be 
zero-free by general principles. However, this seems to be a very challenging 
problem, even for $G=\Sl(2,\R)$. \qed

\def\entries{
\[Ach99 Achab, D., {\it Espace de Hardy pour les quotients
$\Gamma\backslash G$}, Math.\ Z.\ {\bf 230(1)} (1999), 21--45

\[BaOt73 Barth, W., and M.\ Otte, {\it Invariante holomorphe Funktionen auf reduktiven Lie-Gruppen}, 
Math.\ Ann.\ {\bf 201} (1973), 97--112

\[Bo63 Borel, A., {\it Compact Clifford-Klein forms of symmetric spaces}, 
Topology {\bf 2}\linebreak  (1963), 111--122

\[Bo66 ---, {\it Introduction to automorphic forms},
Proc.\ Sympos.\ Pure Math. {\bf 9} (1966), Amer.\ Math.\ Soc., Providence, RI, 199--210

\[Bo69 ---, ``Introduction aux groupes arithm\'etiques'', Hermann, Paris, 1969

\[Fr83 Freitag, E., ``Siegelsche Modulfunktionen'', Grundlehren der Mathematischen Wissenschaften 
{\bf 254}, Springer, 1983

\[GiHu78 Gilligan, B., and A. Huckleberry, {\it On non-compact complex 
nil-manifolds}, Math.\ Ann.\ {\bf 238} (1978), 39--49

\[Hel78 Helgason, S., ``Lie Groups, Differential Geometry and Symmetric Spaces'', Academic Press, London, 
1978

\[HiKr98 Hilgert, J., and B.\ Kr\"otz, {\it The Plancherel Theorem 
for invariant 
Hilbert spaces}, Math.\ Z. {\bf 37(1)} (2001), 31--59 

\[HiKr99 Hilgert, J., and B.\ Kr\"otz, {\it Representation, characters, and spherical functions associated to 
causal symmetric spaces}, J.\ Funct.\ Anal. {\bf 169} (1999), 357--390

\[HiKr01 ---, {\it The Plancherel Theorem 
for invariant Hilbert spaces}, Math.\ Z. {\bf 37(1)} (2001), 31--59 

\[HiNe93 Hilgert, J., and K.-H. Neeb, ``Lie semigroups and their Applications", 
Lecture Notes in Math.\ {\bf 1552}, Springer, 1993

\[H\"o73 H\"ormander, L., ``An introduction to complex analysis in several 
variables'', North-Holland, 1973

\[HoMo79 Howe, R., and C.\  Moore, {\it Asymptotic Properties of Unitary 
Representations}, J.\ Funct.\ Anal.\ {\bf 32} (1979), 72--96 

\[Ko98 Kobayashi, S., ``Hyperbolic complex spaces'', Grundlehren der mathematischen 
Wissenschaften {\bf 318}, Springer 1998

\[Kr98a Kr\"otz, B., {\it On Hardy and Bergman spaces on 
complex Ol'shansk\u\i{} semigroups}, Math.\ Ann.\ {\bf 312} (1998), 13--52

\[Kr98b ---, {\it On the dual of complex Ol'shanski\u{\i}
semigroups}, Math.\ Z. {\bf 237} (2001), 505--529

\[Kr99 ---, {\it The Plancherel Theorem for Biinvariant Hilbert Spaces}, 
Publ.\ RIMS {\bf 35 (1)} (1999),  91--122

\[Kr01 ---, {\it Formal dimension for semisimple symmetric spaces}, Comp. Math.
{\bf 125} (2001), 155--191

\[KN\'O97 Kr\"otz, B., K.-H. Neeb, and G.\ \'Olafsson, {\it  Spherical 
Representations and Mixed Symmetric Spaces}, Represent.\ 
Theory {\bf 1} (1997), 424--461

\[KN\'O01 ---,  {\it  Spherical Functions on Mixed Symmetric Spaces}, 
Represent.\ Theory {\bf 5} (2001), 43--92

\[KO01 Kr\"otz, B., and M.\ Otto, {\it Vanishing properties of analytically continued 
matrix coefficients}, J. Lie Theory, to appear 

\[Lo84 Loeb, J.-J., {\it Fonctions plurisousharmoniques sur un groupe de Lie 
complexe invariantes par une forme r\`eelle}, C. R. Acad. Sci. Paris S\`er. I Math. {\bf 299} (1984), 
no. {\bf 14}, 663--666

\[MaMo60 Matsushima, Y., and A.\ Morimoto, {\it Sur certain espaces 
fibr\'es holomorphe sur une vari\'et\'e de Stein}, 
Bull.\ Soc.\ Math.\ France  {\bf 88}(1960), 137--155

\[Ma97 Mayer, M., {\it Asymptotics of Matrix Coefficients and Closures
of Fourier-Stieltjes Algebras}, J.\ Funct.\ Anal.\ {\bf 143(1)}, (1997), 42--54

\[Ne97 Neeb, K.-H., {\it A general non-linear convexity theorem}, 
Forum Math.\ {\bf 9} (1997), 613--640

\[Ne98 --- , {\it On the complex and convex geometry 
of Ol'shanski\u\i{} semigroups}, Ann.\ Inst.\ Fourier {\bf 48(1)}
(1998), 149--203

\[Ne99a ---, ``Holomorphy and Convexity in Lie Theory,'' 
Expositions in Mathematics {\bf 28}, de Gruyter, 1999

\[Ne99b ---, {\it On the complex geometry of invariant domains in
complexified symmetric spaces}, 
Ann.\ Inst.\ Fourier {\bf 49(1)} (1999), 177--225

\[Ols82 Ol'shanski\u\i, G. I., {\it Invariant cones in Lie algebras, 
Lie semigroups, and the holomorphic discrete series}, Funct. Anal. and Appl. 
{\bf 15}, 275--285 (1982)

\[Ra72 Raghunathan, M.\ S., ``Discrete Subgroups of Lie groups'', 
Ergebnisse der 
Mathematik {\bf 68}, Springer, 1972 

\[Sh92 Shabat, B.\ V., ``Introduction to Complex Analysis, Part II: Functions of Several  
Variables'', Amer. Math. Soc., Providence, Rhode Island, 1992

\[Sta86 Stanton, R. J., {\it Analytic Extension of the holomorphic discrete 
series}, Amer. J. Math. {\bf 108} (1986), 1411--1424

\[WaWo93 Wallach, N., and J. A.  Wolf, {\it  Completeness of Poincar\'e series for automorphic forms associated to
the integrable discrete series}, Representation theory of reductive groups (Park City, Utah, 1982), 265--281, Progr.
Math., 40, Birkhäuser Boston, Boston, Mass., 1983

\[WeWo77 Wells, R. O., and J. A.  Wolf, {\it  Poincar\'e series and automorphic cohomology on flag domains}, 
Ann. of Math. (2) {\bf 105} (1977), no. 3, 397--448

}

{\sectionheadline{\bf References}
\frenchspacing
\entries\par}
\dlastpage
\vfill\eject
\bye